\documentstyle[11pt]{article}

\newtheorem{thm}{Theorem}[section]
\newtheorem{lem}[thm]{Lemma}
\newtheorem{cor}[thm]{Corollary}
\newtheorem{pro}[thm]{Proposition}

\newtheorem{defi}[thm]{Definition}

\newcommand{\gm }{\Gamma }
\newcommand{\lon }{\longrightarrow }
\newcommand{\be }{\begin{eqnarray*}}
\newcommand{\ee }{\end{eqnarray*}}

\newcommand{\st }[1]{\stackrel{#1}{\lon }}

\newcommand{\poiddd }[3]{ (#1\gpd #2, \alpha_{#3}, \beta_{#3})}

\input amssym.def
\input amssym

\newcommand{\pf}{\noindent{\bf Proof.}\ }
\newcommand{\qed}{\begin{flushright} $\Box$\ \ \ \ \ \
                    \end{flushright}}
\newcommand{\complex}{{\Bbb C}}
\newcommand{\reals}{{\Bbb R}}

\newcommand{\frakg}{{\frak g}}

\newcommand{\hstar}{*_{\hbar}}

\newcommand{\half}{\frac{1}{2}}

\newcommand{\F}{{\cal F}}

\newcommand{\cald}{{\cal D}}
\newcommand{\calf}{{\cal F}}

\newcommand{\smalcirc}{\mbox{\tiny{$\circ $}}}

\def\description label#1{\hfil\bf[#1]\hfil}
\newcommand{\secd}[1]{\gm(\wedge^{#1}A^{*})}

\newcommand{\parr}[1]{\frac{\partial  #1}{\partial \lambda_{i}}}
\newcommand{\parrr}[2]{\frac{\partial #1}{\partial #2}}

\newcommand{\rt}{\Lambda}
\newcommand{\rr}{\Lambda^{\#}}

\newcommand{\al}{\alpha}
\newcommand{\lrw}{\longrightarrow}
\newcommand{\ot}{\mbox{$\otimes$}}
\newcommand{\Map}{\longmapsto}

\newcommand{\otr}{\ot_{R}}
\newcommand{\oth}{\ot_{\hbar}}
\newcommand{\otrh}{\ot_{R_{\hbar}}}
\newcommand{\dh}{\cald [[\hbar ]]}
\newcommand{\dhh}{\cald_{\hbar}}
\newcommand{\limh}{\mbox{lim}_{\hbar \mapsto 0}}

\newcommand{\hopf}{(H, R, \alpha , \beta , m , \Delta , \epsilon )}
\newcommand{\eendr}{End_{k}R}
\newcommand{\phia}{\phi_{\alpha}}
\newcommand{\phib}{\phi_{\beta}}
\newcommand{\uqg}{U_{\hbar }\frakg}
\newcommand{\otrf}{\otimes_{R_{\F}}}

\newcommand{\dif}{ (\Delta \otr  id )\F  \F^{12}}
\newcommand{\idf}{(id \otr  \Delta) \F  \F^{23} }
\newcommand{\ohh}{O(\hbar^{3}) }
\newcommand{\oh}{O(\hbar^{2})}
\newcommand{\cdoth}{\cdot_{\hbar}}
\newcommand{\alphah}{\alpha_{\hbar}}
\newcommand{\cdotf}{\cdot_{\F}}

\newcommand{\alt}{\mbox{Alt}}
\newcommand{\fl}{\bar{\Lambda }}
\newcommand{\calfh}{\calf_{\hbar}}
\newcommand{\C}{{\cal M}}
\newcommand{\matroid}{{matrix}}

\renewcommand{\theenumi}{{\rm{(\roman{enumi})}}}

\def\Clin{C^\infty_{\ell in}}

\def\mvf{multiplicative vector field}
\def\mmvf{multiplicative multivector field}
\def\cdo{covariant differential operator}
\def\LAgpd{${\cal LA}$-groupoid}
\def\VBLalgd{${\cal VB}$-Lie algebroid}
\def\VBgpd{${\cal VB}$-groupoid}
\def\CDO{\mathop{\rm CDO}}
\def\End{\mathop{\rm End}}

\def\ST{\ \vert\ }

\def\sdp{\mathbin{\hbox{$\mapstochar\kern-.3333em\times$}}}
\def\pds{\mathbin{\hbox{$\times\kern-.55em\mapstochar\,$}}}

\newcommand{\wed}{\mathbin{\lower1.5pt\hbox{$\scriptstyle{\wedge}$}}}

\let\Tilde=\widetilde
\let\Bar=\overline
\let\Vec=\overrightarrow
\let\ceV=\overleftarrow
\let\Hat=\widehat
\let\isom=\cong
\let\sol=\bullet
\def\chigh{{\raise1.5pt\hbox{$\chi$}}}
\let\phi=\varphi
\def\til0{\Tilde{0}}

\def\dpl{\mathbin{+\hskip-6pt +\hskip4pt}}
\def\dminus{\raise2pt\hbox{\vrule height1pt width 2ex}\hskip3pt}
\def\dtimes{\mathbin{\hbox{\huge.}}}

\def\llangle{\langle\!\langle}
\def\rrangle{\rangle\!\rangle}

\def\pback#1{\mathbin{{{\lower1.2ex\hbox{$\times$}}\atop #1}}}

\def\ddt#1{\left.\frac{d}{dt}#1\right|_0}

\def\brev{\:\breve{\rule{0pt}{8pt}}}
\def\hatt{\Hat{\phantom{X}}}

\def\vlra{\hbox{$\,-\!\!\!-\!\!\!-\!\!\!-\!\!\!-\!\!\!
-\!\!\!-\!\!\!-\!\!\!-\!\!\!-\!\!\!\longrightarrow\,$}}

\def\vleq{\hbox{$\,=\!\!\!=\!\!\!=\!\!\!=\!\!\!=\!\!\!
=\!\!\!=\!\!\!=\!\!\!=\!\!\!=\!\!\!=\!\!\!=\!\!\!=\!\!\!=\,$}}

\def\lrah{\hbox{$\,-\!\!\!-\!\!\!
-\!\!\!-\!\!\!-\!\!\!-\!\!\!-\!\!\!\longrightarrow\,$}}

\def\surj{-\!\!\!-\!\!\!-\!\!\!\gg}

\def\inj{>\!\!\!-\!\!\!-\!\!\!-\!\!\!>}

\def\gpd{\,\lower1pt\hbox{$\longrightarrow$}\hskip-.24in\raise2pt
               \hbox{$\longrightarrow$}\,}

\def\lgpd{\,\lower1pt\hbox{$\vlra$}\hskip-1.02in\raise2pt\hbox{$\vlra$}\,}

\def\llgpd{\,\lower1pt\hbox{$\vvlra$}\hskip-1.3in\raise2pt\hbox{$\vvlra$}\,
}

\def\vgpd{\Bigg\downarrow\!\!\Bigg\downarrow}

\def\vequals{\Bigg\|}

\begin{document}

\title{{\bf Quantum groupoids}}
\author{ PING XU \thanks{ Research partially supported by   NSF
          grants   DMS97-04391 and DMS00-72171.}\\
   Department of Mathematics\\The  Pennsylvania State University \\
University Park, PA 16802, USA\\
          {\sf email: ping@math.psu.edu }}

\date{}

\maketitle

\centerline{{Dedicated to the  memory of  Mosh\'e Flato}}

\begin{abstract}
We introduce a general notion of quantum universal enveloping  algebroids (QUE
algebroids), or quantum groupoids, as a unification of quantum groups and
star-products.  Some  basic properties are studied including   the twist
 construction and the  classical limits. In particular,
we show that a quantum groupoid naturally  gives rise to
a Lie bialgebroid as a classical limit.  Conversely, we formulate
a conjecture on the existence of a quantization for any
Lie bialgebroid, and prove this conjecture for
the special case of  regular triangular Lie bialgebroids.
As an application of this theory, we study the dynamical quantum
groupoid $\cald\ot_{\hbar}\uqg$, which gives an interpretation
of the quantum dynamical  Yang-Baxter  equation in terms of
Hopf algebroids.
\end{abstract}

\section{Introduction}
Poisson  tensors in many aspects  resemble classical triangular r-matrices
in quantum group theory. A  notion unifying both Poisson structures
and Lie bialgebras was introduced in \cite{MX94}, called
     Lie bialgebroids.
The integration theorem for Lie bialgebroids encompasses both Drinfeld's
theorem of integration of Lie bialgebras on the one hand,
and the Karasev-Weinstein theorem of existence  of  local symplectic groupoids
for Poisson manifolds on the other hand \cite{MX98}. Quantization of
Lie bialgebras leads  to quantum groups, while quantizations of Poisson
manifolds are  the so called   star-products.
It is therefore natural to expect that there  exists  some intrinsic
connection between these two quantum objects.
The  purpose  of this paper  is to  fill in this gap
by introducing  the  notion of quantum universal enveloping
 algebroids (QUE algebroids), or quantum groupoids,
as a general framework  unifying   these two concepts.
Part of the results in this paper has been announced in
\cite{Xu1} \cite{Xu2}.

The  general notion of Hopf algebroids was  introduced
by Lu \cite{lu}, where the axioms were
obtained  essentially by translating those of Poisson groupoids
to their quantum counterparts.   The  special
case
where the base algebras are commutative was  studied earlier
    by Maltsiniotis \cite{M}, in a 1992 paper based on the work of
   Deligne on Tannakian categories \cite{Deligne}.   Subsequently,
Vainerman \cite{Vai}  found  a class of examples of Hopf algebroids
arising from a  Hopf algebra action  on an algebra, which
generalizes that introduced by Maltsiniotis. 
Recently,
   Hopf algebroids also appeared in Etingof and Varchenko's work on
   dynamical  quantum groups \cite{EV2} \cite{EV3}.
In this paper, we will  mainly follow Lu's definition,
but some axioms will be modified. One advantage
of our approach  is that the tensor product of
representations becomes an immediate consequence of the definition.
  We refer  the interested
reader to \cite{M1} \cite{NV} for  other
 definitions of quantum groupoids, which are
originated from different motivations and different from
the one we are using here.

As we know, many important examples of Hopf   algebras
    arise as deformations of the universal enveloping algebras of   Lie
algebras.
Given  a Lie algebroid $A$,
its universal enveloping algebra $UA$ (see the definition in Section 2)
    carries   a  natural cocommutative
Hopf algebroid structure.  For example, when $A$ is the
tangent bundle Lie algebroid $TP$, one obtains a  cocommutative
Hopf algebroid structure on   $\cald (P)$, the algebra  of differential
operators on $P$. It is natural to expect that
deformations of $UA$, called quantum universal enveloping algebroids
or quantum groupoids in this paper, would give us some non-trivial
      Hopf algebroids.  This is the  starting point of
    the present  paper.
   Examples include the usual quantum universal enveloping
algebras and
the quantum groupoid $\cald_{\hbar} (P) $ corresponding
to a star-product on a Poisson manifold $P$.

Another important class of quantum groupoids is connected with
the so called quantum
dynamical Yang-Baxter equation, also
known as the Gervais-Neveu-Felder equation \cite{BBB}:
\begin{equation}
\label{eq:dybe}
R^{12}(\lambda )R^{13}(\lambda +\hbar h^{(2)} )R^{23}(\lambda )
=R^{23}(\lambda +\hbar h^{(1)} )R^{13}(\lambda )R^{12}(\lambda +\hbar h^{(3)} ).
\end{equation}
   Here $R(\lambda )$ is a meromorphic function
from $\eta^*$ to $\uqg \ot \uqg$, $\uqg$ is a quasi-triangular
quantum group, and $\eta \subset \frakg$ is an Abelian Lie subalgebra.
This equation  arises
naturally from  various contexts in mathematical physics,
including quantum Liouville theory, quantum
Knizhnik-Zamolodchikov-Bernard equation, and quantum
Caloger-Moser model \cite{ABE} \cite{Babelon} \cite{Felder}.
   One   approach to this equation,
due to Babelon et al. \cite{BBB}, is to use Drinfeld's
theory of quasi-Hopf algebras \cite{dr:quasi}. Consider a meromorphic function
$F: \eta^* \lon \uqg \ot \uqg$ such that $F(\lambda)$ is
invertible for all $\lambda$.  Set
$R(\lambda )=F^{21}(\lambda )^{-1}R F^{12}(\lambda ), $
where $R \in \uqg \ot \uqg$ is the standard universal $R$-matrix for the
quantum group $\uqg$. One can check \cite{BBB} that
$R(\lambda )$ satisfies Equation (\ref{eq:dybe}) if
$F(\lambda )$ is of zero weight,
and satisfies the following shifted cocycle condition:
\begin{equation}
\label{eq:shifted0}
(\Delta_0  \ot  id )F(\lambda )  F^{12} (\lambda +\hbar h^{(3)})
   =  (id \ot  \Delta_{0} )  F (\lambda )F^{23}(\lambda  ),
\end{equation}
where $\Delta_0$ is the coproduct of $\uqg$.
If moreover we assume  that
\begin{equation}
\label{eq:co0}
   (\epsilon_{0} \ot id) F(\lambda )  =  1; \ \
(id \ot \epsilon_{0}) F (\lambda ) =  1,
\end{equation}
where $\epsilon_{0}$ is the counit map,
one can form an elliptic  quantum group, which is a family of
   quasi-Hopf algebras $(\uqg , \Delta_{\lambda})$
   parameterized by $\lambda \in \eta^*$:
$\Delta_{\lambda}=F(\lambda )^{-1}\Delta_{0} F (\lambda )$.
For $\frakg =\frak{sl}_{2}(\complex )$, a solution to Equations
(\ref{eq:shifted0}) and (\ref{eq:co0})
   was obtained by Babelon  \cite{Babelon}
in 1991.
For general simple Lie algebras, solutions were recently
   found  independently by Arnaudon et al. \cite{Arnaudon}  and Jimbo
et al. \cite{Jimbo}  based on the approach of
Fr\o nsdal \cite{F}.
Equivalent solutions are also found by Etingof and
Varchenko \cite{EV3} using intertwining operators.
Recently, using a
 method similar to \cite{Arnaudon} \cite{F}
\cite{Jimbo}, Etingof et al. found  a large class of shifted cocycles
\cite{ESS} 
as quantization of  the classical dynamical $r$-matrices of
 semisimple Lie algebras
in  Schiffmann's classification  list \cite{S}.  On the
other hand, for an arbitrary Lie algebra, a  general  recipe
was obtained very recently by the author 
 for finding the shifted cocycles
quantizing the so called 
 splittable classical triangular dynamical $r$-matrices \cite{Xu4}.

As we will see in Section 7, Equation (\ref{eq:shifted0})
   arises naturally from the  ``twistor" equation
of a quantum groupoid. This leads to  another interpretation
of an elliptic quantum group, namely
   as a  quantum groupoid. Roughly speaking, our construction goes as follows.
Instead of $\uqg$, we start with    the algebra
$H=\cald\ot \uqg $, where $\cald$ denotes the algebra of meromorphic
   differential operators on $\eta^*$.
$H$ is no longer a Hopf algebra. Instead it is a QUE algebroid
   considered as the  Hopf algebroid tensor product of
$\cald $ and $\uqg$.
Then the  shifted  cocycle condition is
shown to be equivalent to the equation  defining a
twistor of this Hopf algebroid.
Using this twistor, we obtain
   a new QUE algebroid $\cald\ot_{\hbar }\uqg$ (or a quantum groupoid).
We note that
$\cald\ot_{\hbar}\uqg$ is co-associative as  a Hopf
algebroid, although $(\uqg, \Delta_{\lambda})$ is
not co-associative. The   construction of
$\cald\ot_{\hbar}\uqg$
is in some sense to restore  the co-associativity  by enlarging the algebra
   $\uqg$ by tensoring the dynamical part $\cald$.
The relation between this quantum groupoid and
     quasi-Hopf algebras $(\uqg , \Delta_{\lambda})$
is, in a certain sense,  similar  to that between
   a fiber bundle and
its fibers.
   We expect that this quantum groupoid will be
useful in understanding elliptic quantum groups, especially
their representation theory \cite{Felder}.
The physical meaning of it,   however,  still
needs to be explored.

This paper is organized as follows: In Section 2, we recall some
basic definitions and results  concerning Lie bialgebroids.
Section 3 is devoted to the definition and basic properties
of Hopf algebroids. In particular, for Hopf algebroids
with anchor, it is proved that  the category of left modules
   is a monoidal category. As a fundamental construction, in Section 4,
we study the twist construction of Hopf algebroids, which
generalizes the usual twist construction of Hopf algebras.
Despite its complexity compared to Hopf algebras,
the fundamentals are analogous   to those of Hopf algebras.
In particular, the monoidal  categories of left modules of
the twisted and untwisted Hopf algebroids are  always
equivalent. Section 5 is devoted to the introduction
of quantum universal enveloping algebroids.
The main part is to show that Lie bialgebroids indeed
appear as the classical limit of QUE algebroids, as
   is expected. However, unlike the quantum group
case, the proof is not trivial and  is in fact rather involved.
On the other hand, the inverse question: the
   quantization problem,
which would   encompass  both  quantization of Lie bialgebras
and deformation quantization of Poisson manifolds as special cases,
remains widely open.  As a very special case,
in Section 6, we show that any regular
triangular Lie bialgebroid is quantizable.
The discussion on quantum groupoids associated to
quantum dynamical $R$-matrices (i.e. solutions
to the  quantum dynamical  Yang-Baxter  equation) occupies Section 7.
The last section, Section 8, consists of an appendix, as well
as a list of open questions.

We would like to mention the recent work by Etingof and
Varchenko \cite{EV2} \cite{EV3},
   where  a different approach to the quantum dynamical
Yang-Baxter equation in the framework of Hopf algebroids was   given.

{\bf Acknowledgments.} The author  would   like to thank  Giuseppe Dito,
Vladimir Drinfeld,  Pavel Etingof,  Masaki Kashiwara,
   Jiang-hua Lu,  Dale Peterson
   and Alan Weinstein for useful discussions and comments.
In addition to the funding source mentioned
in the first footnote, he would also  like to thank RIMS, IHES  and
    Max-Planck Institut for their hospitality and financial support while
part of
this project was being done.

\section{Preliminary on  Lie bialgebroids}
It is well known that the classical objects corresponding to quantum
groups are Lie bialgebras. Therefore, it is not surprising to expect
that the classical  counterparts  of quantum groupoids
are  Lie bialgebroids. However,  unlike Lie bialgebras,  Lie bialgebroids
were introduced and studied before the invention
of quantum groupoids.
In fact, they were used mainly in the study of
Poisson geometry in connection with symplectic
and Poisson groupoids (see \cite{Weinstein} \cite{Xu}).

The purpose of this section is to recall some basic facts  concerning
Lie bialgebroids. We start with recalling some  definitions.

\begin{defi}\label{Lie.Algebroid} A  Lie algebroid
is a vector bundle $A$ over $P$ together with  a Lie
algebra structure on the space $\Gamma(A)$ of smooth sections of $A$,
and  a bundle map $\rho: A \rightarrow TP$ (called the anchor),
   extended to a map between sections of these bundles,  such that

(i) $\rho ([X,Y])=[\rho (X),\  \rho (Y)]$; and

(ii) $[X, fY] = f[X,Y] + (\rho(X) f)Y$\\
   for any smooth sections $X$ and $Y$ of $A$ and any smooth function
$f$ on $P$.
\end{defi}

Among many examples of Lie algebroids are
   usual Lie algebras, the tangent bundle
of a manifold, and an integrable distribution over a manifold
   (see \cite{Mackenzie}).
   Another interesting example is connected with
Poisson manifolds.
Let $P$ be  a  Poisson manifold with Poisson tensor
    $\pi $. Then  $T^{*}P$  carries    a  natural
Lie algebroid structure, called
the cotangent bundle  Lie algebroid of the Poisson manifold $P$ \cite{CDW}.
   The anchor map
   $\pi^{\#}: T^*P \rightarrow TP$ is defined by
\begin{equation}
\label{eq:anchor0}
   \pi^{\#}: ~~T_{p}^{*}P \longrightarrow T_p P:~~
\pi^{\#}(\xi)(\eta )=\pi (\xi ,\eta ), \ \ \ \forall \xi , \eta\in T_{p}^{*}P
\end{equation}
   and the Lie bracket of $1$-forms $\alpha$ and $\beta$ is given by
\begin{eqnarray}
   \label{eq_bracket-on-one-forms}
   [\alpha, \beta  ]
   & = &  L_{{\pi}^{\#}(\alpha)} \beta ~ -
   ~ L_{{\pi}^{\#}(\beta)} \alpha   ~ - ~
d (\pi (\alpha, \beta)).
   \end{eqnarray}

Given a  Lie algebroid $A$, it is known that
$\oplus_{k}\gm (\wedge^{k} A^* )$   admits  a
differential  $d$ that  makes it into a
differential graded algebra \cite{KS-M}. Here, $d: \gm (\wedge^{k} A^* )
\lon \gm (\wedge^{k+1} A^* )$ is  defined by
(\cite{Mackenzie} \cite{MX94}  \cite{WX}):
\begin{eqnarray}
d\omega (X_1,\ldots ,X_{k+1}) & = & \sum_{i=1}^{k+1}  (-1)^{i+1} \rho (X_i)
(\omega (X_1,\hat{\ldots}, X_{k+1})) \nonumber \\
& & \qquad + \sum_{i\ < j}  (-1)^{i+j} \omega ([X_i ,X_j ],X_1 ,
\hat{\ldots}\,\hat{\ldots}\,, X_{k+1}), \label{eq:derivative}
\end{eqnarray}
for $\omega  \in \secd{k}$, $X_i \in \Gamma A,\ 1\leq i\leq k+1$.
Then $d^2 =0$ and one obtains a cochain complex
whose cohomology is called the Lie algebroid  cohomology.
On the other hand,  the Lie bracket on $\gm (A)$
extends naturally to a graded Lie bracket on $\oplus_{k} \gm (\wedge^{k}A)$
called
the Schouten bracket, which,  together with the usual
wedge product $\wedge$, makes  it into a
Gerstenhaber algebra \cite{Xu3}.

As in the case of Lie algebras, associated to any
Lie algebroid, there is an associative algebra
called  the {\em  universal enveloping algebra} of the Lie algebroid $A$
\cite{Huebschmann}, a concept whose definition we now recall.

Let  $A\to P$ be a Lie algebroid  with anchor $\rho$. Then
the $C^{\infty}(P)$-module direct sum
$C^{\infty}(P)\oplus \gm(A)$ is  a Lie algebra over $\reals $  with
the Lie bracket:
$$[f+X, \ g+Y]=(\rho(X ) g -\rho(Y) f)+[X, Y].$$

Let $U=U(C^{\infty}(P)\oplus \gm(A))$ be its  standard
universal enveloping algebra.
For any $f\in C^{\infty}(P)$ and $X\in \gm (A)$, denote
by $f'$ and $X'$ their canonical image in  $U$.
Denote by $I$ the two-sided  ideal  of $U$  generated by all
elements of the form $(fg)'-f'g'$ and $(fX)'-f'X'$.
Define $U(A)=U/I$,   which is
called the universal enveloping algebra of the Lie algebroid $A$.
When $A$ is a Lie algebra, this definition reduces to the
definition of usual  universal enveloping algebras.
On the other hand,  for the tangent
bundle Lie  algebroid  $TP$,   its universal enveloping
   algebra  is $\cald (P)$, the algebra of differential operators over $P$.
In between, if $A=TP\times \frakg$ as the Lie algebroid
direct product, then $U(A)$ is isomorphic to $\cald (P)\ot U\frakg$.
Note that the maps $f \mapsto f'$ and $X \mapsto X'$ considered
   above descend to
linear embedings  $i_{1}: C^{\infty}(P)\mapsto U(A)$, and
$i_{2}: \gm (A)\mapsto U(A)$; the first map $i_1$ is an
algebra morphism. These maps have the following properties:
\begin{equation}
\label{eq.relations}
i_{1}(f)i_{2}(X)=i_{2}(fX), \ \
[i_{2}(X),i_{1}(f)]=i_{1}(\rho (X)f)), \ \
[i_{2}(X),i_{2}(Y)]=i_{2}([X,Y]).
\end{equation}
In fact, $U(A)$ is universal among all  triples $(B, \phi_{1}, \phi_{2})$
having such  properties (see \cite{Huebschmann} for a proof of this
simple  fact).
Sometimes, it is also   useful  to  think  of $UA$    as the  algebra of  left
invariant differential operators on
a local Lie groupoid $G$ which  integrates  the Lie algebroid $A$.

The notion of Lie bialgebroids is a  natural generalization of
that of Lie bialgebras. Roughly speaking, a Lie bialgebroid
is a pair of Lie algebroids ($A$, $A^*$) satisfying
a certain compatibility condition.
    Such a condition, providing
a definition of {\em Lie bialgebroids}, was given in
\cite{MX94}.  We quote here an equivalent formulation from
\cite{K-S:1994}.
\begin{defi}
A {\em Lie bialgebroid} is a dual pair $(A,A^*)$ of vector bundles equipped
with
Lie algebroid structures such that the differential $d_*$ on
$\Gamma(\wedge^*A)$
coming from the structure on $A^*$ is a derivation of the Schouten bracket
on
$\Gamma(\wedge^*A)$.
Equivalently, $d_{*}$ is a derivation for sections of $A$, i.e.,
\begin{equation}
\label{eq:1}
d_{*}[X, Y]=[d_{*}X, Y]+[X, d_{*}Y], \ \ \forall X ,Y \in \gm (A).
\end{equation}
In other words,
   $(\oplus_{k} \Gamma(\wedge^k A), \wedge , [\cdot ,\cdot ], d_{*})$
is a strong differential  Gerstenhaber algebra \cite{Xu3}.
\end{defi}

In fact, a Lie bialgebroid is equivalent to a strong differential
Gerstenhaber algebra structure on $\oplus_{k} \Gamma(\wedge^{k} A)$
(see Proposition 2.3 in \cite{Xu3}).

For a Lie bialgebroid $(A, A^{*})$, the base $P$ inherits a natural
Poisson structure:
\begin{equation}
\label{eq:base}
\{f, \ g\}=<df , \ d_{*}g>, \ \ \forall f, g\in C^{\infty}(P),
\end{equation}
which satisfies  the identity: $[df, \ dg]=d\{f, \ g\}$.

As in the case of Lie bialgebras, a useful method of
constructing Lie bialgebroids is via $r$-matrices.
More precisely,  by an {\em $r$-\matroid}, we mean
a section $\Lambda \in  \gm ( \wedge^2 A ) $ satisfying
\begin{equation}
\label{eq:2}
L_{X}[\rt, \rt ]= [ X, [\rt, \rt]] = 0, \  \ \  \forall X \in \Gamma (A).
\end{equation}

An $r$-\matroid\  $\Lambda$ defines a Lie bialgebroid, where the
    differential  $d_{*}: \Gamma (\wedge^{*} A )
\lon \Gamma (\wedge^{*+1} A)$ is simply given by $d_{*} =[  \cdot  , \  \Lambda ]$.
In this case, the bracket on $\gm (A^* )$ is given by
\begin{equation}
\label{eq*}
[\xi , \eta ]=L_{\rr \xi}\eta -L_{\rr \eta}\xi -d[\Lambda  (\xi , \eta )],
\end{equation}
and the anchor is the composition  $\rho \circ \rr :A^*  \lon TP$, where
    $\Lambda^\#$ denotes     the bundle map $A^{*}\lon A$ defined by
$\rr (\xi) (\eta )=\Lambda  (\xi , \eta ), \forall \xi , \eta \in
   \Gamma (A^{*})$.
   Such a Lie bialgebroid is called a {\em  coboundary  Lie bialgebroid},
in
analogy with the Lie algebra case \cite{LX} \cite{LX1}.
It is called a {\em triangular  Lie bialgebroid} if
$[\Lambda , \Lambda ]=0$. In particular it
is called  a {\em regular triangular  Lie bialgebroid} if
$\Lambda$ is of constant rank.

   When $P$ reduces to a point, i.e., $A$ is a Lie
algebra,  Equation (\ref{eq:2}) is equivalent to that $[\rt , \rt ]$ is
$ad$-invariant, i.e, $\rt$ is an $r$-matrix in the ordinary sense.
On the other hand, when $A$ is the tangent
bundle $TP$ with the  standard Lie algebroid structure,
Equation (\ref{eq:2})  is  equivalent
to that $[\rt , \rt ]=0$, i.e., $\rt$ is a Poisson tensor.

Another interesting class of coboundary Lie
bialgebroids is connected with the so called classical dynamical
$r$-matrices.

Let $\frakg$ be  a Lie algebra over $\reals$ (or $\complex$)
and $\eta \subset \frakg$ an Abelian Lie subalgebra.
A classical dynamical $r$-matrix \cite{ABE}\cite{EV1} is a smooth function (or
meromorphic function in the complex case)
   $r: \eta^{*} \lon \wedge^{2} \frakg$ such that\footnote{Throughout
the paper, we follow the sign convention in \cite{ABE} for  the
definition  of a  classical dynamical $r$-matrix in order to be
consistent with the quantum dynamical Yang-Baxter
equation (\ref{eq:dybe}). This  differs a sign from
the one used in \cite{EV1}.}
\begin{enumerate}
\item $r(\lambda )$  is $\eta$-invariant, i.e., $[h, \ r(\lambda )]=0, \
\forall h\in \eta $;
\item $\alt   (dr)-\half [r, r]$ is constant over $\eta^{*}$
with value in $(\wedge^{3} \frakg)^{\frakg}$.
\end{enumerate}
Here $dr$ is considered as a $\eta\otimes \wedge^{2}\frakg$-valued
function over $\eta^*$ and $\alt$ denotes the
standard skew-symmetrization operator.
In particular,  if $\alt   (dr)-\half [r, r]=0$, it is called
a {\em classical triangular dynamical $r$-matrix}.
The following
is a simple example of a classical dynamical $r$-matrix.\\\\
{\bf Example 2.1} Let $\frakg$ be a simple Lie algebra
with root decomposition $\frakg =\eta  \oplus \sum_{\alpha \in \Delta_{+}}
(\frakg_{\alpha}\oplus \frakg_{-\alpha })$, where
$\eta$ is a Cartan subalgebra, and $\Delta_{+}$ is  the set of
positive roots. Then
\[
r(\lambda) \, = 
-\half  \sum _{\alpha \in \Delta_{+}} \,
\coth (\half  \ll \alpha, \lambda \gg) E_{\alpha} \wedge E_{-\alpha},
\]
is a classical dynamical $r$-matrix,  where
$\ll \, , \, \gg$  is the Killing form of $\frakg$,
the $E_{\alpha}$ and
$E_{-\alpha}$'s are standard
   root vectors, and
$\coth (x)  =  {e^x + e^{-x} \over e^x - e^{-x}}$
is the hyperbolic cotangent function.

A classical dynamical $r$-matrix naturally  gives rise to a Lie bialgebroid.

\begin{pro} \cite{BK-S} \cite{LX1}
Let $r:\eta^{*} \lon \wedge^2 \frakg$ be a classical dynamical $r$-matrix.
Then  $A=T\eta^{*} \times \frakg$, equipped
with the standard Lie algebroid structure, together with
$\Lambda =\sum_{i=1}^{k} (\frac{\partial}{\partial \lambda_{i}}
\wedge  h_{i})+r(\lambda )\in \gm (\wedge^{2}A)$
defines a coboundary Lie bialgebroid. Here $\{h_{1}, h_{2}, \cdots ,h_{k}\}$
is a basis of $\eta$, and    $(\lambda_{1}, \cdots ,\lambda_{k})$ is
the induced  coordinate system on $\eta^{*}$.
\end{pro}

We end this section by recalling the definition of
Hamiltonian operators, which will be needed later on.
   Given a Lie bialgebroid $(A, A^* )$
with associated strong differential
Gerstenhaber algebra $(\oplus_{k} \Gamma(\wedge^k A),
   \wedge , [\cdot ,\cdot ], d_{*})$, one may construct a new Lie
bialgebroid via  a twist.
For that, simply let $\tilde{ d_{*}}= d_{*}+[\cdot , \ H ]$ for 
some  $H\in  \gm (\wedge^{2} A)$. It is easy to check \cite{LWX} that
this still defines a   strong differential Gerstenhaber algebra
(therefore a Lie bialgebroid),   if and only if the following
Maurer-Cartan type equation holds:
\begin{equation}
d_{*}H+\half [H, \ H]=0.
\end{equation}
Such an $H$   is called a Hamiltonian operator of the Lie bialgebroid
$(A, A^* )$.

Finally we note that even though we are mainly dealing with
real vector bundles and real
Lie algebroids in this paper, one may also consider
(as already suggested by the early  example  of classical dynamical $r$-matrices)
complex Lie algebroids and  complex Lie bialgebroids. In that case,
one may have to
    use sheaves of  holomorphic sections etc. instead of global ones.
Most  results in this section  still hold after 
 suitable modifications.

\section{Hopf algebroids}

\begin{defi}
\label{dfn_algebroid}
A  Hopf algebroid $(H, R, \alpha ,
   \beta , m , \Delta , \epsilon)$ consists of the following data:

1) a  {\bf total algebra} $H$ with product $m$, a  base algebra $R$,
    a  {\bf source map}: an algebra homomorphism $\al: ~ R \lrw H$,
and  a {\bf target map}: an algebra anti-homomorphism
$\beta: ~ R \lrw H $
such that the images of $\al$ and $\beta$ commute in $H$, i.e., $\forall
a, b \in R$,
$\al(a) \beta (b) ~ = ~ \beta (b) \al (a)$.
There is then a natural $(R ,  R)$-bimodule structure on $H$ given by
   $a \cdot  h = \al (a)h$ and
$ h \cdot  a = \beta (a) h$.
Thus, we can form the $(R, R)$-bimodule
product $H \otr  H$.  It is easy to see that $H \otr  H$  again admits
an  $(R, R)$-bimodule structure.  This will allow us to form the
   triple product $H \otr H \otr  H$ and etc.

2) a  {\bf co-product}: an $(R, R)$-bimodule map
$\Delta: ~ H \lrw H \otr  H$
with $\Delta(1) = 1 \ot 1$ satisfying the co-associativity:
\begin{equation}
\label{eq:coassociative}
(\Delta \otr  id_{\scriptscriptstyle H} ) \Delta
~ = ~ (id_{\scriptscriptstyle H} \otr  \Delta)
   \Delta: ~ H \longrightarrow H \otr H \otr H;
\end{equation}

3) the product and the co-product  are
{\bf compatible} in the following sense:
\begin{equation}
\label{eq:compatible1}
\Delta (h)(\beta (a)\ot 1-1\ot \alpha (a))=0,  \  \  \mbox{ in }
H\otr H, \ \forall a\in R \mbox{ and }
h\in H, \mbox{ and }
\end{equation}

\begin{equation}
\label{eq:compatible2}
\Delta (h_{1}h_{2})=\Delta (h_{1})\Delta (h_{2}), \ \ \forall h_{1},
h_{2}\in H,
\ \ \ \ \mbox{(see the  remark below)};
\end{equation}

4) a  co-unit map: an (R, R)-bimodule map
$\epsilon: ~~ H \lrw R$
satisfying
$\epsilon(1_{\scriptscriptstyle H}) ~ = ~ 1_{\scriptscriptstyle R} $
(it follows then that $\epsilon \beta   =  \epsilon
\al  =  id_{\scriptscriptstyle R}$) and
\begin{equation}
\label{eq_co-unit}
   (\epsilon \otr id_{\scriptscriptstyle H}) \Delta ~ = ~
(id_{\scriptscriptstyle H} \otr \epsilon) \Delta ~ =
~ id_{\scriptscriptstyle H}: ~ H \lrw H.
\end{equation}
Here  we have used the identification:  $R\otr H\cong H\otr R\cong H$
(note that both  maps on the left hand sides of Equation  (\ref{eq_co-unit})
are well-defined).
\end{defi}
{\bf Remark.} It is clear that any left $H$-module is automatically
an $(R, R)$-bimodule.
Now given  any  left $H$-modules $M_{1}$ and $M_{2}$,
define,
\begin{equation}
\label{eq:module-product}
h\cdot (m_{1}\otr m_{2})=\Delta (h)(m_{1}\ot m_{2}), \ \ \ \
\forall h\in H, \ m_{1}\in M_{1}, \ m_{2}\in M_{2}.
\end{equation}
The
right-hand side is a well-defined element in $M_{1}\otr M_{2}$ due to Equation
(\ref{eq:compatible1}).
In particular, when  taking  $M_{1}=M_{2}=H$, we see that the
right-hand side of
Equation (\ref{eq:compatible2}) makes sense.
In fact,   Equation (\ref{eq:compatible2}) implies
   that $M_{1}\otr M_{2}$ is again
a left $H$-module under the action  defined by
Equation (\ref{eq:module-product}). Left $H$-modules
are also called {\em representations} of the Hopf
algebroid $H$ (as an associative algebra). The
category of representations of $H$ is denoted by Rep$H$.

There is an equivalent version for the compatibility condition
3)  due to   Lu   \cite{lu}:
\begin{pro}
\label{pro:Lu}
    The compatibility condition  3) (Equations
(\ref{eq:compatible1})  and (\ref{eq:compatible2}))
is   equivalent to  that the kernel of the  map
\begin{equation}
\label{eq_psi}
\Psi:~
H \ot H \ot H \lrw H \otr H: ~~  \sum h_1 \ot h_2 \ot h_3
   \Map \sum (\Delta h_1) (h_2 \ot h_3)
\end{equation}
is a  left ideal of $H \ot H^{op} \ot H^{op}$, where $H^{op}$
denotes $H$ with the opposite product.
\end{pro}
\pf  Assume that $Ker \Psi$ is a left ideal. It is
clear that  for any $a\in H$, $1\ot \beta (a)\ot 1-1\ot 1\ot \alpha (a)\in
Ker \Psi$. Hence $h\ot \beta (a)\ot 1-h\ot 1\ot \alpha (a)=
(h\ot 1\ot 1)(1\ot \beta (a)\ot 1-1\ot 1\ot \alpha (a))$ belongs
to  $Ker \Psi$.  That is, $\Delta (h)(\beta (a)\ot 1-1\ot \alpha (a))=0$.
To prove  Equation  (\ref{eq:compatible2}), we assume that
$\Delta h_{2}=\sum_{ij}g_{i}\otr g_{j} $ for some  $g_{i}, g_{j}\in H$.
Then $h_{2}\ot 1 \ot 1-\sum_{ij}1 \ot g_{i}\ot g_{j}\in Ker \Psi$.
This implies that $h_{1}h_{2}\ot 1 \ot 1-\sum_{ij}h_{1}\ot g_{i}\ot g_{j}
\in Ker \Psi$ since $Ker \Psi$ is a left ideal. Hence
$\Delta ( h_{1}h_{2})-\sum_{ij}(\Delta h_{1})(g_{i}\ot g_{j})=0$
in $H\otr H$. I.e.,  $\Delta (h_{1}h_{2})=\Delta (h_{1})\Delta (h_{2})$.

Conversely, assume that Equations (\ref{eq:compatible1}) and
(\ref{eq:compatible2}) hold. Suppose that $\sum h_1 \ot h_2 \ot h_3\in
Ker \Psi$. Then for any $x, y, z\in H$, note that
   in $H \ot H^{op} \ot H^{op}$, we have
$(x\ot y\ot z)\sum (h_1 \ot h_2 \ot h_3 ) =\sum xh_{1}\ot h_{2}y \ot h_{3}z$.
Then
\be
&&\Psi ((x\ot y\ot z)\sum (h_1 \ot h_2 \ot h_3))\\
&=& \sum \Delta (  x h_1) (h_2  y\ot h_3 z)\\
&=&  (\Delta  x) \sum (\Delta h_{1})(h_2 \ot h_3)(y\ot z)\\
&=&0.
\ee
That is, $Ker \Psi$ is a  left ideal. \qed
{\bf Remark.}
(1)  In  \cite{lu}, objects satisfying the above axioms are
called bi-algebroids, while  Hopf algebroids
are referred to those admitting an  antipode.
However,  here   we relax the requirement of    the existence of
   an antipode  for Hopf algebroids,
since   many interesting examples, as shown   below,
often  do not admit an  antipode.

(2) In the classical case, the compatibility
between the Poisson structure and the
groupoid structure implies that the base manifold
is a coisotropic submanifold of the Poisson groupoid \cite{Weinstein}.
    For a Hopf algebroid,  it  would be natural
to expect that the  quantum analogue should hold as well, which means that
   the kernel of $\epsilon$ is a  left ideal of $H$. However,
we are not able to prove this at the moment (note that this extra condition
was
required in the definition in \cite{lu}).\\\\

In most
situations,  Hopf algebroids are equipped with an
additional structure, called an {\em anchor} map.
Let $\hopf$ be a Hopf algebroid (over the ground
   field $k$ of characteristic zero).
By $\eendr$, we denote
the algebra of linear endmorphisms
of $R$ over $k$. It is clear that $\eendr$ is an $(R,R)$-bimodule,
where $R$ acts on it from the left by left multiplication
and  acts from the right by right multiplication. Assume that
$R$ is a left $H$-module and  moreover the representation
$\mu:  H\lon \eendr$  is an $(R,R)$-bimodule map.
For any  $x\in  H$ and $a\in R$, we denote by $x(a)$ the
element $\mu (x) (a)$ in $R$.  Define
\begin{eqnarray}
&&\phia, \phib : (H\otr H)\otimes R\lon H, \nonumber\\
&&\phia (x\otr y\otimes a)=x(a)\cdot y, \ \ \mbox{  and  }
\phib (x\otr  y\otimes a)=x\cdot y(a). \label{eq:phi-alp-beta}
\end{eqnarray}
Here $x, y\in H$, $a\in R$,  and the dot $\cdot$ denotes the
$(R,R)$-bimodule structure on $H$.  Note that $\phia$ and
$\phib$ are well defined since $\mu$ is an $(R, R)$-bimodule map.

\begin{pro}
\label{pro:co-unit}
Under the above assumption, and moreover assume that
\begin{equation}
\phia (\Delta x \otimes a)=x\alpha (a), \ \ \mbox{ and } \ \phib (\Delta x
   \otimes a)=x\beta (a), \ \ \forall x\in H, a\in R.
\end{equation}
Then the map $\tilde{\epsilon}: H\lon R$,  $\tilde{\epsilon} x =x(1_{R})$,
satisfies the co-unit property,
i.e., it
is an $(R,R)$-bimodule map, $\tilde{\epsilon}(1_{H})=1_{R}$,
   and satisfies Equation (\ref{eq_co-unit}).
\end{pro}
\pf That $\tilde{\epsilon}$ is an $(R,R)$-bimodule map
follows from the assumption that the representation
$\mu$ is an $(R,R)$-bimodule map.  It is
clear
that $\tilde{\epsilon}(1_{\scriptscriptstyle H})=1_{\scriptscriptstyle R}$.
To prove Equation (\ref{eq_co-unit}),
assume that $\Delta x=\sum_{i}x_{i}^{(1)}\otr x_{i}^{(2)}$.
Then
\be
   ( \tilde{\epsilon} \otr id_{\scriptscriptstyle H}) \Delta  x &= &
\sum \tilde{\epsilon}(x_{i}^{(1)})\otr x_{i}^{(2)}\\
&= &\sum x_{i}^{(1)} (1_{R})\otr x_{i}^{(2)}\\
&= &\phi_{\alpha}(\Delta x\ot 1_{R})\\
&= &x\alpha (1_{R})\\
&= &x.
\ee
Similarly, we have
$(id_{\scriptscriptstyle H} \otr \tilde{\epsilon}) \Delta x =x$.
\qed

It is thus natural to expect that $\tilde{\epsilon} $ coincides
with the co-unit map.

\begin{defi}
\label{def:anchor}
   Given a Hopf algebroid $\hopf$, an anchor map is a representation
$\mu: H\lon \eendr$, which is  an
   $(R, R)$-bimodule map satisfying
\begin{enumerate}
\item $\phia (\Delta x \otimes a)=x\alpha (a)$ and \ $\phib (\Delta x
   \otimes a)=x\beta (a), \ \ \forall x\in H, a\in R$;
\item $x(1_{R})=\epsilon x, \ \ \forall x\in H$.
\end{enumerate}
\end{defi}
{\bf Remark.} For a Hopf algebra, since $R= k$ and $\eendr \cong k$,
one can simply take  
the counit as the anchor.  In this case, the anchor is in fact equivalent
to the counit map. However, for a Hopf algebroid,
the existence of an anchor map is a stronger assumption than
the existence of a counit. In fact, we can require
axioms 1)-3) in Definition (\ref{dfn_algebroid})
   together with the existence of an anchor map to define a Hopf algebroid
with an anchor.
Then the existence of the counit would  be a
direct consequence according to Proposition (\ref{pro:co-unit}).

Given  any $x=\sum x_{1}\otr x_{2}\cdots \otr x_{n}\in \otr^{n}H$
   and $k$ elements ($k\leq n$)
$a_{i_{1}}, a_{i_{2}}, \cdots , a_{i_{k}}\in R$,
we denote by
   $x(\cdot, \cdots, a_{i_{1}}, \cdots, \cdot, a_{i_{k}}, \cdots, \cdot )$
the element  $\sum x_{1}\otr  \cdots \otr x_{i_{1}}(a_{i_{1}})
\otr \cdots \otr x_{i_{k}}(a_{i_{k}}) \otr \cdots \otr x_{n}$ in
$\otr^{n-k}H$.
The anchor map assumption
guarantees that this is a well-defined element.

\begin{pro}
For any $a , b\in R , x\in H$,
\begin{eqnarray}
& &\alpha (a)(b)=ab, \ \ \beta (a)(b)=ba;  \label{eq:alpab}\\
& &\Delta (x)(a,  b)=x(ab); \label{eq:delxab}\\
& & \epsilon (xy)=x(\epsilon (y)). \label{eq:epxy}
\end{eqnarray}
\end{pro}
\pf  To prove Equation (\ref{eq:alpab}),
   we note that $\alpha (a)(b)=(a\cdot 1_{H})(b)
=a( 1_{H} (b))=ab$, where we used the fact that $\mu$ is
an $(R, R)$-bimodule map. Similarly, we have $\beta (a)(b)=ba$.

For Equation  (\ref{eq:delxab}), we have
\be
\Delta (x) (a,  b) &=&\phi_{\alpha}(\Delta x \otimes a)(b)\\
&=&(x\alpha (a))(b)\\
&=&x(\alpha (a) (b))\\
&=&x(ab).
\ee
Here the last step used Equation (\ref{eq:alpab}).

Finally, using   (ii) in Definition (\ref{def:anchor}), we have
$$\epsilon (xy)=(xy)(1_{R})=x(y(1_{R}))=x(\epsilon (y) ). $$
\qed
{\bf Remark.}  Equation (\ref{eq:alpab}) implies that the
induced $(R, R)$-bimodule structure on $R$, where $R$
is
considered as a  left $H$-module, coincides with
the usual one by (left and right) multiplications. In fact,
this condition is equivalent to requiring that
$\mu$ is an $(R, R)$-bimodule map in Definition (\ref{def:anchor}).
Equation (\ref{eq:delxab})  simply means that
$R\otr R\cong R$ as left $ H$-modules. And the last  equation,
Equation (\ref{eq:epxy}),  amounts to saying that
$\epsilon :H\lon R$ is  a  module map as both $H$ and $R$ are
considered as  left $ H$-modules, where
$H$ acts on $H$ by left multiplication.

The following result  follows immediately from
the definitions.

\begin{thm}
Let $(H, R, \alpha , \beta , m , \Delta , \epsilon)$ be a Hopf algebroid
with anchor $\mu$. Then  the category $Rep H$ of left $H$-modules
equipped with the tensor product $\otr$ as defined by Equation
(\ref{eq:module-product}),
the unit object $(R, \mu)$, and the trivial
associativity isomorphisms: $(M_{1}\otr M_{2})\otr M_{3}\lon
M_{1} \otr (M_{2}\otr M_{3})$ is a monoidal category.
\end{thm}
{\bf Example 3.1}  Let $\cald$ denote  the algebra
   of all differential operators
on a smooth manifold
$P$, and $R$  the algebra  of smooth
functions on $P$. Then  $\cald$ is a Hopf algebroid over  $R$.
Here, $\alpha =\beta $ is  the embedding $R\lon \cald $,
while the coproduct $\Delta : \cald \lon \cald \otr \cald$
is defined  by
\begin{equation}
\Delta (D)(f, g)=D(fg), \ \  \ \forall D \in \cald , \mbox{ and } \ f , \
g\in R.
\end{equation}
Note that  $\cald \otr \cald $  is simply the space of bidifferential
operators. Clearly, $\Delta $ is co-commutative, i.e., $\Delta^{op}=\Delta$.
The usual action of differential operators on $C^{\infty}(P)$ defines
an anchor  $\mu :\cald\lon \eendr$.
In this case,
     the co-unit  $\epsilon : \cald \lon R$  is
the natural projection to its $0th$-order part
of  a differential operator.
It is easy to see that   left
$\cald$-modules are $D$-modules in the usual sense,
and the tensor product is the usual tensor product
of $D$-modules over   $R$.
We note, however, that this Hopf algebroid does not admit
an antipode in any natural sense.  Given a differential operator
$D$,  its antipode, if
it exists, would be the dual operator $D^*$. However, the latter is
a differential operator on $1$-densities, which
does not possess  any canonical identification with   a differential operator
on $R$. \\\\

   The construction   above
   can  be generalized to show that the universal enveloping algebra
$UA$ of a Lie algebroid  $A$ admits a  co-commutative Hopf algebroid
structure.

Again  we take $R=C^{\infty}(P)$, and let
    $\alpha =\beta : R\lon  UA$ be the  natural embedding.
For the  co-product, we set
\be
\Delta (f)&=&f\otr 1, \ \  \forall f\in R;\\
\Delta (X)&=& X\otr 1+1\otr X, \ \ \forall X\in \gm (A).
\ee
This formula  extends  to a co-product $\Delta : UA\lon
UA\otr UA$ by  the compatibility condition:
Equations (\ref{eq:compatible1}) and (\ref{eq:compatible2}).
Alternatively, we  may identify
$UA$ as the    subalgebra of $\cald (G)$  consisting of  left
   invariant differential operators on a (local) Lie groupoid $G$
integrating $A$, and then   restrict the co-product $\Delta_{G}$ on $\cald
(G)$
to this subalgebra. This is well-defined since
$\Delta_G  $ maps  left  invariant differential
operators to  left  invariant bidifferential
operators. Finally, the map $(\mu x)(f)=(\rho x)(f), \ \ \forall x\in UA,
f\in R$ defines  an anchor, and the co-unit map is  then
the projection $\epsilon :UA\lon R$, where $\rho : UA\lon \cald (P)$ denotes
the algebra homomorphism extending the anchor of the Lie algebroid
(denoted by the same symbol $\rho$).

\begin{thm}
$(UA, R, \alpha , \beta , m , \Delta , \epsilon)$ is a co-commutative
Hopf algebroid with anchor $\mu$.
\end{thm}

\section{Twist construction}
As in the  Hopf algebra case,
the twist construction is an important method
of producing   examples of Hopf algebroids. This section is devoted to
the study on this  useful construction. We start with the following

\begin{pro}
Let $\hopf$ be  a Hopf algebroid    with anchor $\mu$,  and
let $\phi_{\alpha}$ and $\phi_{\beta}$ be the maps
defined by Equation (\ref{eq:phi-alp-beta}). Then
for any $x, y, z\in H$ and $a\in R$,
\begin{eqnarray}
\phi_{\alpha }((\Delta x)(y\otr  z)\otimes a)&=&x\alpha (y(a))z;
\label{eq:phia}\\
\phi_{\beta  }((\Delta x)(y\otr z)\otimes a)&=&x\beta (z(a))y. \label{eq:phib}
\end{eqnarray}
\end{pro}
\pf Assume that $\Delta x=\sum_{i}x_{i}^{(1)}\otr x_{i}^{(2)}$. Then
\be
\phi_{\alpha}((\Delta x)(y\otr z)\otimes a)&=&
\phi_{\alpha}(\sum_{i}(x_{i}^{(1)}y\otr x_{i}^{(2)}z)\otimes a)\\
&=&\sum_{i}\alpha ((x_{i}^{(1)}y)(a))x_{i}^{(2)}z\\
&=&\sum_{i}\alpha (x_{i}^{(1)}(y(a)))x_{i}^{(2)}z.
\ee
On the other hand,   using (i) in  Definition \ref{def:anchor},
\be
x\alpha (y(a))z&=&\phi_{\alpha}(\Delta x\otimes y(a))z\\
&=&\phi_{\alpha}(\sum_{i}x_{i}^{(1)}\otr x_{i}^{(2)}\otimes y(a))z\\
&=&\sum_{i}\alpha (x_{i}^{(1)}(y(a))) x_{i}^{(2)}z.
\ee

Hence, $\phi_{\alpha }((\Delta x)(y\otr z)\otimes a)=x\alpha (y(a))z$.
Equation (\ref{eq:phib})  can be proved similarly. \qed

Now  let $\F$ be  an element in $H\otr H$.  Define $\alpha_{\F}, \beta_{\F}:
R\lon H$, respectively, by
\begin{equation}
\alpha_{\F}(a)=\phia (\F \ot a), \ \ \ \beta_{\F}(a)=\phib (\F \ot a),
\ \forall a\in R.
\end{equation}
And for any $a, b\in R$, set
\begin{equation}
a*_{\F}b=  \alpha_{\F}(a)(b).
\end{equation}
More explicitly, if $\F=\sum_{i}x_{i}\otr y_{i}$ for $x_{i}, y_{i}\in H$,
then  $ \forall a, b\in R$,
\begin{eqnarray}
&&\alpha_{\F}(a)=\sum_{i} x_{i}(a) \cdot y_{i}=\sum_{i} \alpha
(x_{i}(a))y_{i},
\label{eq:alpha}\\
&& \beta_{\F}(a)= \sum_{i} x_{i} \cdot y_{i}(a)=\sum_{i}\beta (
y_{i}(a))x_{i},
\ \ \ \mbox{  and } \label{eq:beta}\\
&& a*_{\F}b=\sum_{i} x_{i}(a) y_{i}(b) \label{eq:aFb}.
\end{eqnarray}

\begin{pro}
\label{pro:twist}
Assume that $\F\in H\otr H$   satisfies:
\begin{eqnarray}
&& \dif =\idf \ \ \mbox{ in } \ \  H \otr H \otr H; \ \mbox{ and
}  \label{eq:cocycle}\\
&& (\epsilon \otr id) \F ~ = ~ 1_{H}; \ \
(id \otr \epsilon) \F ~ = ~ 1_{H}.
\label{eq:co-unit}
\end{eqnarray}
   Here $\F^{12}=\F\otimes 1\in (H\otr H)\ot H$,   $\F^{23}=1\otimes \F
\in H\ot (H\otr H)$, and in Equation (\ref{eq:co-unit})
we have used the identification:  $R\otr H\cong H\otr R\cong H$
(note that both  maps on the left hand sides of Equation (\ref{eq:co-unit})
are well-defined). Then
\begin{enumerate}
\item $(R, *_{\F})$ is an associative algebra, and
$1_{R}*_{\F}a=a*_{\F}1_{R}=a, \ \  \forall a \in R$.
\item $\alpha_{\F} :R_{\F}\lon H$ is an algebra homomorphism,
and $\beta_{\F} :R_{\F}\lon H$ is
an algebra anti-homomorphism.
Here   $R_{\F}$ stands for the algebra $(R, *_{\F})$.
\item $(\alpha_{\F} a)(\beta_{\F}b)=(\beta_{\F}b)(\alpha_{\F} a), \ \ \forall
a, b\in R \ .$
\end{enumerate}
\end{pro}
\pf As a first step, we prove that for any $a, b\in R$,
\begin{eqnarray}
\alpha_{\F}(a*_{\F}b)&=&(\alpha_{\F}  a ) (\alpha_{\F} b), \label{eq:alpf}\\
\beta_{\F}(a*_{\F}b)&=&(\beta_{\F}  b ) (\beta_{\F} a). \label{eq:betaf}
\end{eqnarray}

Assume that $\F=\sum_{i}x_{i}\otr y_{i}$ for $x_{i}, y_{i}\in H$. Then
\begin{eqnarray}
\dif &=&\sum_{ij}\Delta x_{i}(x_{j}\otr y_{j})\otr y_{i},\label{eq:dif}\\
\idf &=&\sum_{ij}x_{i}\otr \Delta y_{i}(x_{j}\otr y_{j}). \label{eq:idf}
\end{eqnarray}

Thus
\be
[\dif ](a, b, \cdot )
&=&\sum_{ij} \Delta x_{i} (x_{j}\otimes y_{j})(a, b)\otr y_{i}\\
&=&\sum_{ij}x_{i}(x_{j}(a)y_{j}(b))\otr y_{i}\\
&=&\sum_{ij}\alpha [x_{i}(x_{j}(a)y_{j}(b))]y_{i}\\
&=&\alpha_{\F}(\sum_{j}x_{j}(a)y_{j}(b))\\
&=&\alpha_{\F}(a*_{\F}b),
\ee
where the second equality used Equation (\ref{eq:delxab}).

On the other hand,

\be
[\idf ] (a, b, \cdot )
&=&\sum_{ij} x_{i}(a)\otr \phi_{\alpha }(\Delta y_{i}(x_{j}\otr y_{j}) \otimes
b)\\
&=&\sum_{ij} x_{i}(a)\otr y_{i}\alpha (x_{j}(b))y_{j}\\
&=&\sum_{ij} \alpha ( x_{i}(a))y_{i}\alpha (x_{j}(b))y_{j}\\
&=&(\alpha_{\F}a )( \alpha_{\F}b).
\ee
Thus Equation (\ref{eq:alpf}) follows from Equation (\ref{eq:cocycle}).
The equation $\beta_{\F}(a*_{\F}b)=(\beta_{\F}  b )( \beta_{\F} a)$
can be proved similarly.

Now   for  any $a, b, c\in R$,
   $[(\alpha_{\F}  a )( \alpha_{\F} b)](c)=(\alpha_{\F}  a)( (\alpha_{\F}
b)(c))
=a*_{\F}  (b*_{\F} c)$.
On the other hand, $\alpha_{\F}(a*_{\F}b)(c)=(a*_{\F}b )*_{\F}c$.
The associativity of $R_{\F}$  thus follows from Equation (\ref{eq:alpf}).

Finally, we have
$\alpha_{\F}  (1_{R})=\sum_{i}x_{i}(1_{R})\cdot y_{i}
=\sum_{i} \epsilon (x_{i} )\cdot y_{i}
= (\epsilon \otimes_{R}id)\F
= 1_{H}$. Similarly, $\beta_{\F} (1_{R})=1_{H}$.
It thus follows that
$1_{R}*_{\F}a=\alpha_{\F}(1_{R})(a)=1_{H}(a)=a$.
Similarly, $a*_{\F}1_{R}=a$.

For the last statement,  a similar computation leads to
\be
&&[ {\dif} ](a, \cdot , b)=(\beta_{\F}b ) ( \alpha_{\F}a), \ \ \mbox{ and}\\
&&[{\idf} ] (a, \cdot , b)=(\alpha_{\F}a ) ( \beta_{\F}b ).
\ee
Thus (iii) follows immediately.  This concludes the proof. \qed

\begin{pro}
   Under the same
hypotheses as in Proposition \ref{pro:twist}, we have
\begin{equation}
\label{eq:f2}
\F   (\beta_{\F} (a)\otimes 1-1\ot \alpha_{\F} (a))=0 \ \mbox{ in } H\otr
H,
   \ \ \forall a\in R.
\end{equation}
   \end{pro}
\pf
   \be
&&\F   (\beta_{\F} (a)\otimes 1-1\ot \alpha_{\F} (a))\\
&=&(\sum_{i}x_{i}\otr y_{i})(\sum_{j}\beta (y_{j}(a))x_{j}\otimes 1-
\sum_{j}1\otimes \alpha (x_{j}(a))y_{j})\\
&=&\sum_{ij}[x_{i}   \beta (y_{j}(a))x_{j}\otr y_{i}-
x_{i}\otr y_{i} \alpha (x_{j}(a))y_{j}].
\ee

Now using Equations (\ref{eq:dif})-(\ref{eq:idf}) and
 (\ref{eq:phia})-(\ref{eq:phib}), we obtain
\be
&&[\dif ] (\cdot , a , \cdot )\\
&=&\sum_{ij}\phi_{\beta}(\Delta x_{i}(x_{j}\otr y_{j})\otimes a)\otr y_{i}\\
&=&\sum_{ij}x_{i}\beta (y_{j}(a))x_{j}\otr y_{i},
\ee
and
\be
&& [\idf ](\cdot , a , \cdot )\\
&=&\sum_{ij}x_{i}\otr \phi_{\alpha}( \Delta y_{i}(x_{j}\otr y_{j})\otimes a)\\
&=&\sum_{ij}x_{i}\otr y_{i}\alpha (x_{j}(a))y_{j}.
\ee
Thus the conclusion follows  immediately   from Equation
(\ref{eq:cocycle}).\qed

As an immediate consequence, we have

\begin{cor}
\label{cor:tensor}
Let $M_{1}$ and $M_{2}$ be any  left $H$-modules. Then
\begin{eqnarray}
\F^{\#} : M_{1}  \otrf M_{2}&\lon  &M_{1}\otr M_{2} \nonumber \\
(m_{1}\otrf m_{2})&\lon &\F \cdot (m_{1}\ot m_{2}),\ \  m_{1}\in M_{1}, \
\mbox{ and } \ m_{2}\in M_{2}, \label{eq:fn}
\end{eqnarray}
is  a  well defined linear map.
\end{cor}

Note that $M_{1}\otr M_{2}$ is automatically
an $(R, R)$-bimodule since both $M_1$ and
$M_2$ are $(R, R)$-bimodules. Similarly, $M_{1}  \otrf M_{2}$
is an $(R_{\F}, R_{\F})$-bimodule.
Besides,
   $M_{1}\otr M_{2}$ is also a left $H$-module.
The next lemma indicates how  these module structures are related.

\begin{lem}
\label{lem:fam}
For any $a\in R$ and $m\in M_{1}  \otrf M_{2}$,
\begin{eqnarray}
\F^{\#}(a\cdotf m)&=&\alpha_{\F}(a)\cdot \F^{\#}(m);\\
\F^{\#}(m\cdotf a )&=&\beta_{\F}(a)\cdot \F^{\#}(m),
\end{eqnarray}
where the dot on the
 right-hand side
means the left $H$-action
on $M_{1}\otr M_{2}$, and the dot $\cdotf$  on the left-hand side refers
to
both the left and right  $R_{\F}$-actions on $M_{1}  \otrf M_{2}$.
\end{lem}
\pf  For simplicity, let us assume that $m=m_{1}\otrf m_{2}$
for $m_{1}\in M_{1}$ and $m_{2}\in M_{2}$. Then
\be
&&\F^{\#}(a\cdotf (m_{1}\otrf m_{2}))\\
&=&\F^{\#} ((\alpha_{\F}(a)  m_{1})\otrf m_{2})\\
&=&\F \cdot (\sum_{i}\alpha (x_{i}(a))y_{i}m_{1}\otrf m_{2})\\
&=&\sum_{ij} x_{j}\alpha (x_{i}(a))y_{i}m_{1}\otr y_{j}m_{2}\\
&=&\dif (a, \cdot , \cdot)\cdot (m_{1}\ot m_{2}) .
\ee

On the other hand,

\be
&&\alpha_{\F}(a) \cdot \F^{\#}(m_{1}\otrf m_{2})\\
&=& \alpha_{\F}(a) \cdot (\sum_{j} x_{j} m_{1}\otr  y_{j}m_{2})\\
&=&\sum_{j}\Delta (\alpha_{\F}(a))(x_{j} m_{1}\otr y_{j}m_{2})\\
&=&\sum_{ij}\Delta (\alpha (x_{i}(a))y_{i})( x_{j} m_{1}\otr  y_{j}m_{2})\\
&=&\sum_{ij} (x_{i}(a) \cdot \Delta y_{i})( x_{j} m_{1}\otr y_{j}m_{2})\\
&=&\sum_{ij}[(x_{i}(a) \otr \Delta y_{i}( x_{j} \otr y_{j})](m_{1}
\otimes m_{2})\\
&=&\idf (a, \cdot , \cdot )\cdot (m_{1} \otimes m_{2}).
\ee

The conclusion thus follows again  from Equation (\ref{eq:cocycle}). \qed

We     say that  $\F $ is {\em invertible} if
$\F^{\#}$  defined by Equation (\ref{eq:fn}) is a vector space isomorphism
for any  left $H$-modules $M_{1}$ and $M_{2}$. In this case,
   in particular we can take $M_{1}=M_{2}=H$  so that  we have
an isomorphism
\begin{equation}
\label{eq:f}
\F^{\#} : H\otrf H \lon H \otr H .
\end{equation}

An immediate consequence of Lemma \ref{lem:fam}  is the following

\begin{cor}
\label{cor:f}
If $\F$ is invertible, then for any $a\in R_{\F}$ and
   $n\in  M_{1}\otr M_{2}$,
\be
&&\F^{\#-1} (\alpha_{\F}(a) \cdot n)=a\cdotf \F^{\#-1} (n);\\
&&\F^{\#-1} (\beta_{\F}(a) \cdot n)=\F^{\#-1} (n) \cdotf a.
\ee
\end{cor}

\begin{defi}
An element $\F \in H\otr H$ is called a twistor if
it is invertible and  satisfies  Equations  (\ref{eq:cocycle}) and
(\ref{eq:co-unit}).
\end{defi}

Now assume that $\F$ is a twistor.
Define a new coproduct  $\Delta_{\F}: H \lon H\otrf H $ by
\begin{equation}
\label{eq:coproduct}
\Delta_{\F}=\F^{-1}\Delta \F ,
\end{equation}
where Equation (\ref{eq:coproduct}) means that
$\Delta_{\F}(x)=\F^{\# -1}
(\Delta (x)\F )$, $\forall x \in H$.
In what follows, we will prove that $\Delta_{\F}$ is
indeed a Hopf algebroid co-product.

\begin{lem}
\label{lem:delf}
For any $x\in H$ and $a\in R_{\F}$,
\begin{eqnarray}
&&\Delta (a\cdot_{\F} x)=\alpha_{\F}(a)\cdot \Delta x; \label{eq:dax}\\
&&\Delta (x\cdot_{\F} a)=\beta_{\F}(a) \cdot \Delta x, \label{eq:dxa}
\end{eqnarray}
where $\cdot_{\F}$ refers to
both the  left and the
  right $R_{\F}$-actions on $H$.
\end{lem}
\pf  We have
\be
&&\Delta (a\cdot_{\F} x)\\
&=&\Delta ( \alpha_{\F}(a)x)\\
&=&\Delta ( \alpha_{\F}(a))\Delta x\\
&=& \alpha_{\F}(a)\cdot \Delta x.
\ee
Equation (\ref{eq:dxa}) can be proved similarly.
\qed

\begin{pro}
$\Delta_{\F}: H\lon H\otrf H $  is an $(R_{\F}, R_{\F})$-bimodule map.
\end{pro}
\pf   For any $a\in R_{\F}$ and $x\in H$, using Lemma \ref{lem:delf}
and Corollary \ref{cor:f}, we have
\be
&&\Delta_{\F}(a\cdot_{\F}x)\\
&=&\F^{\# -1}[\Delta (a\cdot_{\F}x )\F]\\
&=&\F^{\# -1}[\alpha_{\F}(a)\cdot \Delta x \F]\\
&=&a\cdotf  \F^{\# -1} (\Delta x \F)\\
&=&a \cdot_{\F} \Delta_{\F} (x).
\ee
Similarly, we can show that
   $\Delta_{\F}(x\cdot_{\F}a )=\Delta_{\F} x \cdot_{\F}a $.\qed

\begin{pro}
\label{pro:com}
The comultiplication $\Delta_{\F}: H\lon H\otrf H $  is compatible with the
multiplication in $H$.
\end{pro}
\pf Consider the maps
\begin{equation}
\label{eq_phi}
\Psi:~
H \ot H \ot H \lrw H \otr H: ~~  \sum h_1 \ot h_2 \ot h_3
   \Map \sum (\Delta h_1) (h_2 \ot h_3),
\end{equation}
and

\begin{equation}
\label{eq_phif}
\Psi_{\F}:~
H \ot H \ot H \lrw H \otrf H: ~~  \sum h_1 \ot h_2 \ot h_3
   \Map \sum (\Delta_{\F} h_1) (h_2 \ot h_3).
\end{equation}

We first prove that
\begin{equation}
\label{eq:psif}
\F^{\#}\smalcirc \Psi_{\F}=\Psi \smalcirc \F_{23}^{\#}.
\end{equation}
\be
(\F^{\#}\smalcirc \Psi_{\F})(\sum h_1 \ot h_2 \ot h_3)&=&
\sum \F^{\#}( (\Delta_{\F} h_1) (h_2 \ot h_3))\\
&=&\sum \F^{\#}[{\F^{\#}}^{-1} ((\Delta h_{1})\F )( h_2 \ot h_3)]\\
&=&\sum \Delta h_{1} \F (h_2 \ot h_3)\\
&=&(\Psi \smalcirc \F_{23}^{\#})(\sum h_1 \ot h_2 \ot h_3).
\ee

Equation (\ref{eq:psif}) implies that
$Ker \Psi_{\F}=Ker \Psi$. To see this,
note that $x\in Ker \Psi_{\F}$, i.e., $ \Psi_{\F}(x)=0$, is
equivalent to $(\F^{\#}\smalcirc \Psi_{\F}) (x)=0$, which is
equivalent to $(\Psi \smalcirc \F_{23}^{\#})(x)=0$, or
$\F_{23}^{\#}(x)\in  Ker \Psi$. Since $Ker \Psi$
is a left ideal in $H \ot H^{op} \ot H^{op}$, the latter is
equivalent to
the fact that   $x\in  Ker\Psi$.
The final conclusion thus follows from Proposition \ref{pro:Lu}. \qed

Proposition \ref{pro:com} implies that $M_{1}\otrf M_{2}$
is again
a left $H$-module for any  left $H$-modules $M_{1},  M_{2}$.
Moreover, it is easy to see that $\F^{\#}: M_{1}\otrf M_{2}\lon  M_{1}\otr
M_{2}$
is an isomorphism of left  $H$-modules. Now we are ready to prove that
$\Delta_{\F}$ is coassociative.

\begin{pro}
$\Delta_{\F}: H\lon H\otrf H $  is coassociative.
\end{pro}
\pf Assume that $M_{1}, M_{2}, M_{3}$ are any left
$H$-modules.  It suffices to prove that the natural
identification  $\phi_{\F}:
(M_{1} \otrf M_{2})\otrf  M_{3} \lon M_{1} \otrf (M_{2}\otrf  M_{3})$
is an isomorphism of  left $H$-modules.

Consider the following diagram:
\be
   \begin{array}{lllll}
(M_{1} \otrf M_{2})\otrf  M_{3} & \st{(\F^{12})^{\#} }
&(M_{1} \otr M_{2})\otrf
   M_{3}&\st{[(\Delta \otr id)\F ]^{\#}}&(M_{1} \otr M_{2})\otr  M_{3} \\\\
\phi_{\F} \downarrow  & & & &\downarrow  \phi \\\\
M_{1} \otrf (M_{2}\otrf  M_{3}) & \st{(\F^{23})^{\#} }
&M_{1} \otrf (M_{2}\otr  M_{3})&\st{[(id\otr \Delta) \F ]^{\#} }
&M_{1} \otr (M_{2} \otr  M_{3})
   \end{array} \ee
   Equation (\ref{eq:cocycle}) implies  that the above diagram commutes.
Since all the other maps involved  in the diagram above  are
isomorphisms of left  $H$-modules, $\phi_{\F}$ is an $H$-module isomorphism
as well. This concludes the proof. \qed

By now, we have actually proved all the Hopf algebroid axioms for
$(H, R_{\F}, \alpha_{\F} , \beta_{\F}, m, \Delta_{\F}, \epsilon)$ except
for the condition on counit $\epsilon$. Instead of proving
this last condition directly, here we show that $\mu$ is still an anchor
after the twist, and therefore Axiom 4) in
Definition \ref{dfn_algebroid} would be a consequence
according to the remark  following Definition
\ref{def:anchor}.  Note that $R_{\F}$ can still
be considered as a   left $H$-module under the  representation
$\mu : H\lon End R_{\F}$ (here only the underlying vector space structure
on $R_{\F}$ is involved). We prove that $\mu$ still  satisfies
the anchor axioms.

\begin{lem}
\label{lem:phiF}
For any $x, y\in H$ and $a\in R$,
\begin{eqnarray}
\phi_{\alpha }^{\F}(  (x   \otrf  y)\otimes a)&=&
   \phi_{\alpha }( \F^{\#} (x   \otrf  y)\otimes a);\\
\phi_{\beta }^{\F}(  (x   \otrf  y)\otimes a)&=&
   \phi_{\beta  }( \F^{\#} (x   \otrf  y)\otimes a).
\end{eqnarray}
\end{lem}
\pf
\be
\phi_{\alpha }^{\F}(  (x   \otrf  y)\otimes a)&=&x(a)\cdot_{\F} y\\
&=&\alpha_{\F}(x(a))y\\
&=&\sum_{i} \alpha (x_{i}(x(a)))y_{i}y\\
&=&\sum_{i} \alpha (x_{i}x(a))y_{i}y.
\ee
On the other hand,
\be
&&\phi_{\alpha }( \F^{\#} (x   \otrf  y)\otimes a)\\
&=&\phi_{\alpha}(\sum_{i} x_{i}x\otr y_{i}y\ot a)\\
&=&\sum_{i} \alpha ( x_{i}x (a))y_{i}y.
\ee
Hence, $\phi_{\alpha }^{\F}(  (x   \otrf  y)\otimes a)=
   \phi_{\alpha }( \F^{\#} (x   \otrf  y)\otimes a)$. Similarly,
one can prove that $\phi_{\beta }^{\F}(  (x   \otrf  y)\otimes a)=
   \phi_{\beta  }( \F^{\#} (x   \otrf  y)\otimes a)$.
\qed

\begin{pro}
The map $\mu : H\lon End R_{\F}$ satisfies the anchor
axioms in Definition \ref{def:anchor} for
   $(H, R_{\F}, \alpha_{\F} , \beta_{\F}, m, \Delta_{\F},
\epsilon)$.
\end{pro}
\pf  First we need to show that $\mu$ is an
$(R_{\F}, R_{\F})$-bimodule map. This can be
checked easily since
$(a\cdot_{\F}x)(b)=(\alpha_{\F}(a)x)(b) =\alpha_{\F}(a) (x(b))=
a*_{\F}x(b)$.  Similarly, $(x\cdot_{\F}a)(b)=x(b)*_{\F}a$.

Axiom (ii) in Definition
\ref{def:anchor}  holds automatically since
$\mu$ is
an anchor for $\hopf$. Now according to Lemma \ref{lem:phiF}
\be
&&\phi_{\alpha }^{\F} (\Delta_{\F}x\ot a)\\
&=&\phi_{\alpha }(\Delta x \F\ot a)\\
&=&\sum_{i}\phi_{\alpha }(\Delta x ( x_{i}\otr y_{i})\ot a)\\
&=&\sum_{i}x\alpha (x_{i}(a))y_{i}\\
&=&x \alpha_{\F}(a).
\ee
Here the second from the last equality used Equation (\ref{eq:phia}).
Similarly,
$\phi_{\beta }^{\F} (\Delta_{\F}x\ot a)=x\beta_{\F}(a)$. This concludes
the proof. \qed

In summary, we have proved
\begin{thm}
\label{thm:twist}
   Assume that $\hopf$ is a Hopf algebroid with anchor $\mu$, and
$\F\in H\otr H $ a twistor. Then
$(H, R_{\F}, \alpha_{\F} , \beta_{\F}, m, \Delta_{\F}, \epsilon)$
is a Hopf algebroid, which still admits $\mu$ as
an anchor. Moreover,
its  corresponding monoidal category of left
$H$-modules is equivalent to that of  $\hopf$.
\end{thm}

We say that $(H, R_{\F}, \alpha_{\F} , \beta_{\F}, m, \Delta_{\F},
\epsilon)$
is obtained from $\hopf$ by twisting via $\F$.

The following theorem  generalizes   a  standard result in
Hopf algebras \cite{Drinfeld1}.
\begin{thm}
   If $\F_{1}\in H\otr H$ is a  twistor
for the Hopf algebroid $H$,
   and $\F_{2}\in
H\ot_{R_{\F_{1}}}H$ a twistor for the twisted
Hopf algebroid $H_{\F_{1}}$, then
the Hopf algebroid obtained by twisting $H$  via $\F_{1}$ then
via $\F_{2}$ is equivalent to that obtained by twisting
via $\F_{1} \F_{2}$. Here $\F_{1} \F_{2}\in H\otr H$
is understood as $\F_{1}^{\#}( \F_{2} )$, where
$\F_{1}^{\#} : H\ot_{R_{\F_{1}}}H\lon H\otr H$ is the map as
defined  in  Equation (\ref{eq:f}).
\end{thm}
\pf
   Clearly, $\F=\F_{1} \F_{2}=\F_{1}^{\#}( \F_{2} )$ is a well defined
element in $H\otr H$.  We only need to verify
that $\F$ is still a twistor. The rest of the
theorem follows from  a routine verification.
 For this purpose, it suffices to show that
$\F$ satisfies both Equation (\ref{eq:cocycle})  and 
Equation  (\ref{eq:co-unit}).
To check that,   in fact we may think  of $\F_1$ and $\F_2$
as elements in $H\ot H$ by taking  some  representatives.   Then
\be
&& \dif  \\
&=& (\Delta \otr  id )\F_{1}
(\Delta \otr  id)\F_{2}
\F_{1}^{12}\F_{2}^{12}\\
&=&(\Delta \otr  id )\F_{1}  \F_{1}^{12}
   [(\F_{1}^{12})^{-1} (\Delta \otr  id )\F_{2} \F_{1}^{12}] \F_{2}^{12}\\
&=&[(\Delta \otr  id )\F_{1} \F_{1}^{12}]
[(\Delta_{\F_1} \ot_{R_{\F_{1}}}  id)\F_{2}
     \F_{2}^{12}]\\
&=&[  (id\otr  \Delta) \F_{1}    \F_{1}^{23}]
[(id \ot_{R_{\F_{1}}} \Delta_{\F_1}  )\F_{2}
    \F_{2}^{23}]\\
&=&(id \otr  \Delta) \F_{1}  (id \otr \Delta  )\F_{2} \F_{1}^{23}
\F_{2}^{23}\\
&=&\idf .
\ee

To prove Equation (\ref{eq:co-unit}),
   assume that $\F_{1} =\sum_{i}x_{i}^{(1)}\otr y_{i}^{(1)}$,
and $\F_{2} =\sum_{i}x_{i}^{(2)}\ot_{R_{\F_{1}}} y_{i}^{(2)}$.
Then $\F_{1} \F_{2}=\sum_{ij} x_{i}^{(1)} x_{j}^{(2)}  \otr y_{i}^{(1)}
   y_{j}^{(2)}$.
And
\be
   &&(\epsilon \otr id) \F \\
&=&\sum_{ij}\epsilon (x_{i}^{(1)} x_{j}^{(2)})\otr  y_{i}^{(1)}  y_{j}^{(2)}\\
&=&\sum_{ij}\epsilon (x_{i}^{(1)} x_{j}^{(2)})\cdot
y_{i}^{(1)}  y_{j}^{(2)} \ \ \ (\mbox{using Equation (\ref{eq:epxy})})\\
&=&\sum_{ij}x_{i}^{(1)} (\epsilon x_{j}^{(2)})\cdot y_{i}^{(1)}  y_{j}^{(2)}\\
&=&\sum_{ij} \phi_{\alpha }(x_{i}^{(1)}\otr y_{i}^{(1)}\ot
\epsilon  x_{j}^{(2)}) y_{j}^{(2)}\\
&=&\sum_{j} \phi_{\alpha } (\F_{1}\ot \epsilon  x_{j}^{(2)}) y_{j}^{(2)}\\
&=&\sum_{j}\alpha_{\F_1}(\epsilon  x_{j}^{(2)})y_{j}^{(2)}\\
&=&(\epsilon \ot_{R_{\F_{1}}} id)\F_{2}\\
& = & 1_{H}.
\ee
Similarly, we prove that $(id \otr \epsilon ) \F=1_{H}$.
\qed

We end this section by the following:\\\\
{\bf Example 4.1} Let $P$ be a smooth manifold,   $\cald$ the
algebra  of differential operators on $P$, and $R=C^{\infty}(P)$. Let $\dh$
denote the space of formal power series in $\hbar$ with coefficients
in $\cald$. The Hopf algebroid structure on $\cald$
naturally extends to
   a Hopf algebroid structure on $\dh$ over the base  algebra $R[[\hbar ]]$,
   which   admits a natural anchor map.

Let $\F =1\otr 1+\hbar B_{1}+\cdots \in \cald\otr \cald [[\hbar ]]
(\cong \dh\ot_{R[[\hbar ]]}\dh )$
be a  formal power series
   of bidifferential operators.
It is easy to see that  $\F$ is a twistor iff the multiplication
on $ R[[\hbar ]]$ defined by:
\begin{equation}
\label{eq:star}
f\hstar g=\F (f, g), \ \ \forall f, g\in R[[\hbar ]]
\end{equation}
   is associative with identity  being the constant function $1$, i.e.,
$\hstar $ is a star product on $P$. In this case,
the bracket  $\{f, g\}=B_{1}(f, g)-B_{1}(g, f), \ \ \ \forall
f , \ g\in C^{\infty}(P)$,  defines  a Poisson structure on $P$,  and
$f\hstar g=\F (f, g)$
is simply  a deformation quantization of  this
    Poisson structure \cite{BEFFL}.

The twisted Hopf algebroid can be easily described. Here
   $\dhh=\dh$ is   equipped with the usual multiplication,
$R_{\hbar}=R[[\hbar ]]$ is   the   $*$-product defined by
Equation (\ref{eq:star}), $\alpha_{\hbar} :R_{\hbar}\lon \dhh$
and $\beta_{\hbar} :R_{\hbar}\lon \dhh$ are given,
respectively,  by
$$\alpha_{\hbar} (f)g=f\hstar g, \ \ \  \beta_{\hbar} (f)g=g\hstar f, \ \
\forall f, g\in R. $$
The co-product $\Delta_{\hbar}: \dhh\lon \dhh\ot_{R_{\hbar}}\dhh$ is
$$\Delta_{\hbar}=\F^{-1}\Delta \F,$$
and the  co-unit $\epsilon$ remains the same, i.e., the projection
$\dhh\lon R_{\hbar}$. This twisted Hopf algebroid
$(\dhh , R_{\hbar}, \alpha_{\hbar} , \beta_{\hbar} , m, \Delta_{\hbar},
\epsilon )$
is called the  quantum groupoid associated to the  
 star product $*_{\hbar}$ \cite{Xu1}.




\section{Quantum groupoids and their classical limits}

The main  purpose of this section is to introduce
quantum universal enveloping algebroids (QUE algebroids), also
called quantum groupoids in the paper, as a
deformation of the standard Hopf algebroid $UA$.

\begin{defi}
\label{def:deformation}
A deformation of a Hopf algebroid  $(H, R, \alpha , \beta , m , \Delta ,
\epsilon)$
over a field $k$ is a topological Hopf algebroid
$(H_{\hbar}, R_{\hbar}, \alpha_{\hbar} , \beta_{\hbar} , m_{\hbar} ,
\Delta_{\hbar} , \epsilon_{\hbar})$ over the ring $k[[\hbar ]]$ of
formal power series in $\hbar $ such that
\begin{enumerate}
\item $H_{\hbar }$ is isomorphic to $H[[\hbar ]]$ as $k[[\hbar ]]$ module with
identity  $1_{H}$, and
$R_{\hbar }$ is isomorphic to $R[[\hbar ]]$ as $k[[\hbar ]]$ module with
identity $1_{R}$;
\item  $\alpha_{\hbar}=\alpha (\mbox{mod } \hbar ),\ 
\beta_{\hbar}=\beta (\mbox{mod } \hbar ), \ \ m_{\hbar}=m  (\mbox{mod }
\hbar ),\ 
\epsilon_{\hbar }=\epsilon  (\mbox{mod } \hbar )$;
\item $\Delta_{\hbar}  =\Delta ( \mbox{mod } \hbar ) $.
\end{enumerate}
In this case,  we  simply say that the quotient
      $H_{\hbar}/\hbar H_{\hbar}$ is
isomorphic to $H$ as a   Hopf algebroid.
\end{defi}

Here the meaning of (i) and (ii) is quite clear. However, for
Condition (iii), we need the following simple fact:
\begin{lem}
\label{lem:limit}
Under the
hypotheses (i) and (ii) as in Definition  \ref{def:deformation},
$H_{\hbar}\otrh H_{\hbar} /\hbar (H_{\hbar}\otrh H_{\hbar})$
is isomorphic to $H\otr H$ as a   $k$-module.
\end{lem}
\pf  Define $\tau :H_{\hbar}\ot H_{\hbar}\lon H\otr H$ by
$$(\sum_{i}x_{i}\hbar^{i})\ot (\sum_{i}y_{i}\hbar^{i})  \lon
x_{0}\otr y_{0}.$$
For any $a\in R$ and $x, y\in H$, since
\be
&&(\beta_{\hbar}a\ot 1-1\ot \alpha_{\hbar} a)(x\ot y)\\
&=&(\beta a)x\ot y-x\ot (\alpha a)  y +O(\hbar ),\\
\ee
then $\tau [(\beta_{\hbar}a\ot 1-1\ot \alpha_{\hbar} a)(x\ot y)]=0$.
In other words, $\tau $ descends to a well defined map
from $H_{\hbar}\otrh H_{\hbar}$ to $ H\otr H$.
It is easy to see that $\tau $ is surjective and
   $Ker \tau =\hbar (H_{\hbar}\otrh H_{\hbar}) $. The conclusion thus
follows immediately. \qed

   By abuse of notation, we  still use    $\tau $ to denote the
induced map $ H_{\hbar}\otrh H_{\hbar}\lon H\otr H$.
   We shall also use the
notation $\hbar \mapsto  0$ to denote this map whenever  the
meaning  is clear from the context.
Then, Condition (iii) means that $\limh \Delta_{\hbar} (x)=\Delta (x)$
for any $x\in  H$.

\begin{defi}
A quantum  universal enveloping  algebroid (or
QUE algebroid),
also called a quantum groupoid,  is a deformation of
the  standard  Hopf algebroid  $(UA,  R, \alpha ,
   \beta , m , \Delta , \epsilon)$  of   a Lie algebroid $A$.
\end{defi}

Let $U_{\hbar}A=UA [[\hbar ]]$ and  $R_{\hbar}=R[[\hbar ]]$.
Assume that
   $(U_{\hbar}A,  R_{\hbar},
   \alpha_{\hbar}, \beta_{\hbar}, m_{\hbar}, \Delta_{\hbar},
\epsilon_{\hbar})$
is  a quantum groupoid. Then  $R_{\hbar}$ defines a
star product on $P$ so that  the equation
$$\{f , g\}=\limh \frac{1}{\hbar} (f*_{\hbar}g-g*_{\hbar}f), \ \ \
\forall f, g \in R $$
defines a Poisson structure on the base space  $P$.

Now define
\be
\delta f&=&\limh \frac{1}{\hbar} ( \alpha_{\hbar}f -\beta_{\hbar}f)
\in UA , \ \ \ \forall f\in R, \\
\Delta^{1}X&=&\limh \frac{1}{\hbar} (\Delta_{\hbar}X-(1\otrh X+X\otrh 1))
\in UA \otr UA,
\ \ \ \forall X\in \gm (A),  \ \mbox{ and }\\
\delta X&=&\Delta^{1}X-(\Delta^{1}X)_{21} \in UA \otr UA  .
\ee
Here for $T= \sum x\otr y\in UA \otr UA$, $T_{21} =\sum y\otr x$.
For the convenience of notations, we introduce
$$\alt T=T-T_{21}, \ \ \ \ \forall T\in  UA \otr UA$$
   so that $\delta X =\alt \Delta^{1}X$.

Below we will use $*_{\hbar}$ to denote both the
 multiplication
  in $U_{\hbar}A$ and that in
$R_{\hbar}$ provided there is no confusion.
For any $f, g\in R$, $x, y\in UA$,
    write
\be
\alpha_{\hbar} f&=&f+\hbar \alpha_{1}f +\hbar^{2}  \alpha_{2}f+
O(\hbar^{3});\\
   \beta_{\hbar} f&=&f+\hbar  \beta_{1}f + \hbar^{2}  \beta_{2}f+
O(\hbar^{3});\\
f*_{\hbar}g&=&fg +\hbar B_{1}(f, g)  +O(\hbar^{2});\\
x*_{\hbar}y &=&xy  +\hbar m_{1}(x, y) +O(\hbar^{2}),
\ee
where    $\alpha_{1}f,  \beta_{1}f,
\alpha_{2}f,  \beta_{2}f $ and $m_{1}(x, y)$ are  elements in $UA$. Hence,
\be
\{f, g\}&=&B_{1}(f, g)-B_{1}(g, f), \ \ \mbox{ and}\\
\delta f&=&\alpha_{1}f -\beta_{1}f.
\ee

\begin{lem}
\label{lem:alpha1}
For any $f, g\in R$,
\begin{enumerate}
\item $ \alpha_{1}(fg) =g(\alpha_{1}f)+f(\alpha_{1}g)
+[\alpha_{1}f, g]+m_{1}(f, g) -B_{1}(f, g);$
\item $ \beta_{1}(fg) =g(\beta_{1}f)+f(\beta_{1}g)+[\beta_{1}f, g]+m_{1}(f, g)
-B_{1}(g, f);$
\item $ [\alpha_{1}f, g]-[\beta_{1}g, f]=m_{1}(g, f)-m_{1}(f, g). $
\end{enumerate}
\end{lem}
\pf From the identity
$\alpha_{\hbar}(f*_{\hbar}g )=\alpha_{\hbar}f *_{\hbar} \alpha_{\hbar}g$,
it follows, by considering  the  $\hbar^1$-terms, that
$$\alpha_{1}(fg)+B_{1}(f, g)=m_{1}(f, g)+(\alpha_{1}f)g+f(\alpha_{1}g). $$
   Thus (i)   follows immediately. And   (ii)  can be proved similarly.

On the other hand, we know,  from the
definition of Hopf algebroids, that
$(\alpha_{\hbar}f)*_{\hbar}(\beta_{\hbar}g)=(\beta_{\hbar}g)*_{\hbar}
(\alpha_{\hbar}f). $
By considering the $\hbar^1$-terms,  we obtain
$$(\alpha_{1}f)g+f(\beta_{1}g)+m_{1}(f, g)=g(\alpha_{1}f)
+(\beta_{1}g)f+m_{1}(g,  f). $$
This proves (iii).
\qed

\begin{cor}
For any $f, g\in R$,
\begin{enumerate}
\item $\delta (fg)=f\delta g+g\delta f $;
\item $[\delta f, g] =\{f , g\}$.
\end{enumerate}
\end{cor}
\pf By symmetrizing   the third identity in
   Lemma \ref{lem:alpha1}, we obtain that
$$[\alpha_{1}f-\beta_{1}f, g]-[\beta_{1}g-\alpha_{1}g, f] =0, $$
i.e., $$[\delta f, g]=-[\delta g, f].$$

Now subtracting Equation (ii) from  Equation  (i) in Lemma
\ref{lem:alpha1}, one  obtains that
$$\delta (fg)=g\delta f +  f\delta g +  [\delta f, g]-\{f, g\}.$$
I.e.,
$$\delta (fg)-(g\delta f + f\delta g )= [\delta f, g]-\{f, g\}.$$
Note that  the
left-hand side
of this equation is symmetric with respect to $f$ and $g$,
   whereas the
right-hand side is skew-symmetric,  so both sides  must  vanish. The
conclusion thus follows immediately. \qed

\begin{lem}
\label{lem:3.7}
For any $f\in R$,
\begin{enumerate}
\item $\Delta_{\hbar}f=f\otrh 1+\hbar (\alpha_{1}f\otrh
1-\Delta_{\hbar}\alpha_{1}f)
+\hbar^{2}(\alpha_{2}f\otrh 1-\Delta_{\hbar}\alpha_{ 2}f) +\ohh$;
\item $\Delta_{\hbar}f=1\otrh f+\hbar (1\otrh \beta_{1}f-\Delta_{\hbar}
\beta_{ 1}f) +\hbar^{2}( 1\otrh \beta_{2}f-\Delta_{\hbar}\beta_{ 2}f) +\ohh$;
\item $f\otrh 1-1 \otrh f=\hbar (1\otrh \alpha_{1}f-\beta_{1}f
\otrh 1)+\hbar^{2}(1\otrh \alpha_{2}f-\beta_{2}f \otrh 1) +\ohh$.
\end{enumerate}
\end{lem}
\pf Since $\Delta_{\hbar}: U_{\hbar}A\lon  U_{\hbar}A\otrh  U_{\hbar}A$ is
an $(R_{\hbar}, R_{\hbar})$-bimodule map, it follows
that
$$\Delta_{\hbar}(f\cdoth 1)=f\cdoth \Delta_{\hbar}1. $$
Here, as well as in the sequel, $\cdoth$ denotes
both the left and the right $R_{\hbar}$-actions on $U_{\hbar}A$,
   and  on   its  appropriate tensor powers.

Now $$f\cdoth 1=\alpha_{\hbar}f=f+\hbar\alpha_{1}f+\hbar^{2}\alpha_{2}f
+\ohh, $$
while
$$f\cdoth \Delta_{\hbar}1=f\cdoth (1\otrh 1)=\alpha_{\hbar}f\otrh 1
=(f+\hbar\alpha_{1}f+\hbar^{2}\alpha_{2}f)\otrh 1+\ohh .$$
Thus it follows that
   $$\Delta_{\hbar}f=f\otrh 1+\hbar (\alpha_{1}f\otrh 1-\Delta_{\hbar}\alpha_{
1}f)
+\hbar^{2}(\alpha_{2}f\otrh 1-\Delta_{\hbar}\alpha_{ 2}f) +\ohh .$$
Similarly, one  can prove (ii).

Finally, since $1\otrh \alpha_{\hbar} f=\beta_{\hbar} f\otrh 1$,
we  have
$$1\otrh (f+\hbar\alpha_{1}f+\hbar^{2}\alpha_{2}f +\ohh )
=(f+\hbar\beta_{1}f+\hbar^{2}\beta_{2}f+\ohh )\otrh 1. $$
This  implies (iii). \qed

\begin{cor}
\label{cor:deltaf}
For any $f, g\in R$,
\begin{enumerate}
\item $\Delta (\delta f)=\delta f\otr 1+1\otr \delta f$;
\item $\Delta^{1}(\delta f)=\nu_{2}f\otr 1+1\otr \nu_{2}f-\Delta \nu_{2}f$,
where $\nu_{2}f=\alpha_{2}f -\beta_{2}f$.
\end{enumerate}
\end{cor}
\pf Combining the three identities in Lemma \ref{lem:3.7} ((i)-(ii)+(iii)),
   we obtain that
\be
&&\delta f\otrh 1+1\otrh \delta f-\Delta_{\hbar}\delta f\\
&& \ \ +\hbar [\alpha_{2}f\otrh 1+1\otrh \alpha_{2}f-\beta_{2}f\otrh 1-
1\otrh \beta_{2}f -\Delta_{\hbar}(\alpha_{2}f -\beta_{2}f)]+\oh=0.
\ee

I.e.,
\begin{equation}
\label{eq:tem}
\delta f\otrh 1+1\otrh \delta f-\Delta_{\hbar}\delta f
+\hbar [\nu_{2}f\otrh 1+1\otrh \nu_{2}f -\Delta_{\hbar} \nu_{2}f]+\oh =0.
\end{equation}

By letting $\hbar  \mapsto 0$, this implies that
   $\delta f\otr 1+1\otr \delta f-\Delta(\delta f ) =0$.
This concludes the proof of (i).

Now writing   $\Delta_{\hbar} \delta f=\delta f\otrh 1+1\otrh \delta
f+\hbar \Delta_{\hbar}^{1} (\delta f)$,
and substituting  it back into Equation (\ref{eq:tem}),
we obtain that
$$\Delta_{\hbar}^{1} (\delta f)=
\nu_{2}f\otrh 1+1\otrh \nu_{2}f -\Delta_{\hbar} \nu_{2}f+O(\hbar ). $$
(ii) thus follows immediately by letting $\hbar  \mapsto 0$.\qed

An immediate consequence is

\begin{cor}
For any $f\in  R$, $\delta f\in \gm (A)$ and $\delta^{2}f=0$.
\end{cor}
\pf  From Corollary \ref{cor:deltaf} (i), it follows that
$\delta f$ is primitive, i.e., $\delta f\in \gm (A)$.
According to Corollary  \ref{cor:deltaf} (ii),
   $\Delta^{1}(\delta f)$ is symmetric. Therefore,
$\delta^{2}f=\delta (\delta f)$, being the skew-symmetric
part of $\Delta^{1}(\delta f)$, equals  zero. \qed

\begin{lem}
\label{lem:delta}
For any $X\in \gm (A)$,
\begin{equation}
\label{eq:delta2}
\Delta^{1} X \otr 1+(\Delta\otr id)\Delta^{1} X
=1 \otr \Delta^{1} X +(id \otr \Delta) \Delta^{1} X.
\end{equation}
\end{lem}
\pf For any $X\in \gm (A)$, denote
$$\Delta^{1}_{\hbar}X=\frac{1}{\hbar }[\Delta_{\hbar}X-(1\otrh X+X\otrh
1)]. $$
Thus $\Delta^{1} X=\limh  \Delta^{1}_{\hbar}X$ and
   $$\Delta_{\hbar} X=1\otrh X + X\otrh 1+\hbar \Delta_{\hbar}^{1} X. $$
Then
\be
(\Delta_{\hbar} \otrh id)\Delta_{\hbar} X&=&
1\otrh 1\otrh X+ \Delta_{\hbar} X \otrh 1 +\hbar
(\Delta_{\hbar} \otrh id)\Delta_{\hbar}^{1} X\\
&=& 1\otrh 1\otrh X+ 1\otrh X\otrh 1 +X\otrh 1\otrh 1 +\hbar
[\Delta_{\hbar}^{1} X\otrh 1 +(\Delta_{\hbar}\otrh id )\Delta_{\hbar}^{1}X)];
\ee
and
\be
   (id \otrh  \Delta_{\hbar}) \Delta_{\hbar}X&=&
1\otrh \Delta_{\hbar} X +X\otrh 1\otrh 1+\hbar (id \otrh \Delta_{\hbar})
\Delta_{\hbar}^{1} X\\
&=&1\otrh 1 \otrh X+1\otrh X\otrh 1 +X\otrh 1\otrh 1 +\hbar
[1\otrh \Delta_{\hbar}^{1} X +(id \otrh \Delta_{\hbar}) \Delta_{\hbar}^{1} X].
\ee

It thus follows that

\begin{equation}
\label{eq:delta1}
\Delta_{\hbar}^{1} X \otrh 1+(\Delta_{\hbar}\otrh id )\Delta_{\hbar}^{1} X
=1 \otrh \Delta_{\hbar}^{1} X +(id \otrh \Delta_{\hbar}) \Delta_{\hbar}^{1} X.
\end{equation}
The conclusion thus follows immediately by letting $\hbar \mapsto  0$. \qed

According to Proposition \ref{pro:appendix}
in the Appendix, we immediately have the following

\begin{cor}
For any $X\in \gm (A)$, $\delta X\in \gm (\wedge^{2}A)$.
\end{cor}

\begin{lem}
For any $f\in  R$ and $X\in \gm (A)$,
$$\delta (fX)=f\delta X+\delta f\wedge X .$$
\end{lem}
\pf For any $X\in \gm (A)$,  again we let
$$\Delta^{1}_{\hbar}X=\frac{1}{\hbar }[\Delta_{\hbar}X-(1\otrh X+X\otrh
1)]. $$

Now
$$f\cdoth X=\alphah f *_{\hbar}X =
fX+\hbar [(\alpha_{1}f)X+m_{1}(f, X)]+O(\hbar^{2}).$$

Hence
\begin{eqnarray}
&&\Delta_{\hbar}(f\cdoth X) \nonumber \\
&=&\Delta_{\hbar}(fX) +\hbar [\Delta_{\hbar}
((\alpha_{1}f)X)+\Delta_{\hbar}m_{1}(f, X)]+O(\hbar^{2})  \nonumber \\
&=&1\otrh fX+fX\otrh 1+\hbar [\Delta^{1}_{\hbar}(fX)
+\Delta_{\hbar}
((\alpha_{1}f)X)+\Delta_{\hbar}m_{1}(f, X)]+\oh  \label{eq:a}.
\end{eqnarray}

On the other hand,
\begin{eqnarray}
&& f\cdoth \Delta_{\hbar} X \nonumber\\
&=& f\cdoth (X\otrh 1+1\otrh X+\hbar \Delta^{1}_{\hbar}X) \nonumber\\
&=&[fX+\hbar ((\alpha_{1}f)X+m_{1}(f, X))]\otrh 1+
(f+\hbar \alpha_{1}f)\otrh X+\hbar f \cdoth \Delta^{1}_{\hbar}X+\oh .
   \label{eq:b}
\end{eqnarray}

  From Equations (\ref{eq:a}) and (\ref{eq:b}), it  follows that

\begin{eqnarray}
1\otrh fX-f \otrh X&=&
\hbar [(\alpha_{1}f)X\otrh 1+m_{1}(f, X)\otrh 1+\alpha_{1}f\otrh X+
f \cdoth\Delta^{1}_{\hbar}X \nonumber \\
&&-\Delta_{\hbar}^{1}(fX)-  \Delta_{\hbar}
((\alpha_{1}f)X)-\Delta_{\hbar}m_{1}(f, X)]+\oh . \label{eq:1fx}
\end{eqnarray}

  From the identity $1 \otrh f \cdoth X=(1\cdoth f) \otrh X$,
it follows that
$$1\otrh (fX+\hbar ((\alpha_{1}f)X+m_{1}(f, X)))=
(f+\hbar \beta_{1}f)\otrh X +O(\hbar^{2}). $$

That is,

\begin{equation}
\label{eq:1fx1}
1\otrh fX- f \otrh  X=\hbar [ \beta_{1}f\otrh X- 1\otrh (\alpha_{1}f)X
-1\otrh m_{1}(f, X)] +\oh .
\end{equation}

By comparing Equations (\ref{eq:1fx}) and (\ref{eq:1fx1}), one obtains that
\be
\Delta_{\hbar}^{1}(fX)
&=&f \cdoth\Delta_{\hbar}^{1}X +(\alpha_{1}f-\beta_{1}f)
\otrh X+(\alpha_{1}f)X\otrh 1 +1\otrh (\alpha_{1}f)X
-\Delta_{\hbar}((\alpha_{1}f)X)\\
&&\ \ + m_{1}(f, X)\otrh 1+1\otrh m_{1}(f, X)-\Delta_{\hbar}m_{1}(f,
X)+O(\hbar ).
\ee

Taking  the limit by letting $\hbar \mapsto 0$, we obtain that
\be
\Delta^{1}(fX)& =&f\Delta^{1}X + \delta f
\otr X+(\alpha_{1}f)X\otr 1 +1\otr (\alpha_{1}f)X
-\Delta((\alpha_{1}f)X)\\
&& \ \ +m_{1}(f, X)\otr 1+1\otr m_{1}(f, X)-\Delta m_{1}(f, X).
\ee

The conclusion thus follows by taking the
skew-symmetrization.  \qed

In summary, we have proved the following

\begin{pro}
\label{thm:delta0}
For any $f, g  \in R$ and $X\in \gm (A)$,
\begin{enumerate}
\item $\delta f\in \gm (A)$ and $\delta X\in \gm (\wedge^{2}A)$;
\item $\delta (fg)=f\delta g+g\delta f $;
\item $\delta (fX)=f\delta X+\delta f\wedge X $;
\item $[\delta f , g]=\{f , g\}$;
\item $\delta^2 f=0$.
\end{enumerate}
\end{pro}

Properties (i)-(iii) above   allow us to extend $\delta$ to a
well-defined degree $1$  derivation
$\delta :\gm (\wedge^{*}A)\lon \gm (\wedge^{*+1}A)$.
Below we will show that $(\oplus \gm (\wedge^{*}A) , \wedge , [\cdot ,
\cdot ], \delta )$ is
a strong differential  Gerstenhaber algebra. For this purpose,
it suffices to show that $\delta$ is a derivation
with respect to $[\cdot , \cdot ]$, and $\delta^2 =0$.
We will prove these facts in two separate propositions
below.

\begin{pro}
\label{thm:delta}
For any $X, Y \in \gm (A)$,
\begin{eqnarray}
&&\delta [X, Y]=[\delta X, Y]+[X,  \delta Y ].
\label{eq:bialgebroid}
\end{eqnarray}
\end{pro}
\pf $ \forall X, Y\in \gm (A)$,
$$\Delta_{\hbar}(X*_{\hbar}Y)=\Delta_{\hbar} X *_{\hbar} \Delta_{\hbar}Y
=(1\otrh X+X\otrh 1+\hbar \Delta^{1}_{\hbar}X)*_{\hbar}
(1\otrh Y+Y\otrh 1+\hbar \Delta^{1}_{\hbar}Y). $$
It thus follows that
\begin{eqnarray}
\Delta_{\hbar} [X, Y]_{\hbar}&=&1\otrh [X, Y]_{\hbar}+[X, Y]_{\hbar}\otrh 1
+\hbar [(\Delta^{1}_{\hbar}X) (1\ot Y+Y\ot 1)+(1\ot X+X\ot
1)\Delta^{1}_{\hbar}Y\nonumber \\
&&\ \ \ -(\Delta^{1}_{\hbar}Y) (1\ot X+X\ot 1)-(1\ot Y+Y\ot
1)\Delta^{1}_{\hbar}X]+O(\hbar^{2}). \label{eq:xyh}
\end{eqnarray}
Here   $[X, Y]_{\hbar}=X*_{\hbar}Y-Y*_{\hbar}X$. Then
$[X, Y]_{\hbar}=[X, Y]+\hbar l_{1}(X, Y)+O(\hbar^{2})$, where
$l_{1}(X, Y)=m_{1}(X, Y)-m_{1}(Y, X)$. Hence
\begin{eqnarray}
\Delta_{\hbar} [X, Y]_{\hbar}&=&
\Delta_{\hbar} [X, Y]+\hbar \Delta_{\hbar}l_{1}(X, Y)+O(\hbar^{2}) \nonumber\\
&=&1\otrh [X, Y]+[X, Y] \otrh 1+\hbar (\Delta_{\hbar}^{1}[X, Y]+
   \Delta_{\hbar} l_{1}(X, Y) )+O(\hbar^{2}) . \label{eq:xyh1}
\end{eqnarray}

Comparing  Equation (\ref{eq:xyh1}) with
    Equation  (\ref{eq:xyh}), we obtain that
\be
\Delta_{\hbar}^{1}  [X, Y]&=&-\Delta_{\hbar}l_{1}(X, Y)+
1\otrh l_{1}(X, Y)+l_{1}(X, Y)\otrh 1\\
&&+(\Delta_{\hbar}^{1} X)(1\ot Y+Y\ot 1)+(1\ot X+X\ot 1)\Delta_{\hbar}^{1}Y\\
&&-(\Delta_{\hbar}^{1} Y)(1\ot X+X\ot 1)-(1\ot Y+Y\ot 1)
\Delta_{\hbar}^{1}X+O(\hbar ).
\ee

This implies, by letting $\hbar\mapsto 0$, that
\be
\Delta^{1}  [X, Y]&=&-\Delta l_{1}(X, Y)+
1\otr l_{1}(X, Y)+l_{1}(X, Y)\otr 1\\
&&+(\Delta^{1} X)(1\ot Y+Y\ot 1)+(1\ot X+X\ot 1)\Delta^{1}Y\\
&&-(\Delta^{1} Y)(1\ot X+X\ot 1)-(1\ot Y+Y\ot 1)\Delta^{1}X.
\ee

Equation (\ref{eq:bialgebroid}) thus follows immediately
by taking the skew-symmetrization.  \qed

\begin{pro}
\label{pro:delta2}
For any $X\in \gm (A)$,
$$\delta^{2} X=0  .$$
\end{pro}
\pf  Let
\begin{equation}
J_{\hbar} =(\Delta_{\hbar}\otrh id )\Delta_{\hbar}^{1} X
-(id \otrh \Delta_{\hbar}) \Delta_{\hbar}^{1} X-1 \otrh \Delta_{\hbar}^{1} X +
\Delta_{\hbar}^{1} X \otrh 1.
\end{equation}
  From  Equation (\ref{eq:delta1}), we know that   $J_{\hbar} =0$.
Let $\{e_{i}\in UA\}$  ($e_{0}=1$) be a local basis of $UA$ over the left
module $R$. Assume that $\delta X=\sum Y_{i}\wedge Z_{i}$ with
$Y_{i}, Z_{i}\in \gm (A)$, and
$$\Delta^{1}X=\sum f^{ij}e_{i}\otr e_{j}+\sum (Y_{i}\otr Z_{i}-Z_{i}\otr
Y_{i}), $$
where $f^{ij}\in R$ are symmetric: $f^{ij}=f^{ji}$.
We may also  assume that  $\Delta e_{i}=\sum g_{i}^{kl}e_{k}\otr e_{l}$ with
$g_{i}^{kl}=g_{i}^{lk}\in R$ since $\Delta$ is co-commutative.
Let us  write
\begin{eqnarray}
\Delta^{1}_{\hbar}X&=&\sum f^{ij}\cdoth e_{i}\otrh e_{j}+\sum (Y_{i}
\otrh Z_{i}-Z_{i}\otrh Y_{i})
+\hbar \Delta_{\hbar}^{2}X; \ \ \ \mbox{ and } \label{eq:d1h}\\
\Delta_{\hbar} e_{i}&=&\sum g_{i}^{kl}\cdoth e_{k}\otrh e_{l}
+\hbar \Delta_{\hbar}^{1}e_{i}, \label{eq:dhe}
\end{eqnarray}
for some  $\Delta_{\hbar}^{2}X$ and $\Delta_{\hbar}^{1}e_{i}\in
U_{\hbar}A\otrh
U_{\hbar}A$.

Then
\be
(id \otrh \Delta_{\hbar}) \Delta_{\hbar}^{1} X&=&
\sum (\alpha_{\hbar}f^{ij}*_{\hbar}\beta_{\hbar} g_{j}^{kl}*_{\hbar}e_{i})
\otrh e_{k}\otrh e_{l}+\sum (Y_{i}\otrh Z_{i}\otrh 1
+Y_{i}\otrh 1\otrh  Z_{i}\\
&&- Z_{i} \otrh Y_{i} \otrh 1-Z_{i}\otrh 1 \otrh Y_{i})
+\hbar [\sum f^{ij}\cdoth  e_{i}  \otrh \Delta_{\hbar}^{1} e_{j}\\
&&+\sum Y_{i}\otrh \Delta_{\hbar}^{1}Z_{i}-\sum  Z_{i}\otrh
\Delta_{\hbar}^{1}Y_
{i}
+( id \otrh  \Delta_{\hbar})\Delta_{\hbar}^{2} X],
\ee

and
\be
1 \otrh \Delta_{\hbar}^{1} X &=&
\sum \beta_{\hbar}f^{ij} \otrh e_{i} \otrh e_{j}+\sum (1\otrh Y_{i}\otrh Z_{i}
- 1\otrh Z_{i}\otrh Y_{i})+\hbar (1 \otrh \Delta_{\hbar}^{2}X).
\ee

Similarly, we may write
\begin{eqnarray}
\Delta^{1}_{\hbar}X&=&\sum e_{i}\otrh e_{j} \cdoth f^{ij}+\sum (Y_{i}
\otrh Z_{i}-Z_{i}\otrh Y_{i})
+\hbar \tilde{\Delta}_{\hbar}^{2}X; \ \ \ \mbox{ and } \label{eq:dlh-tilde} \\
\Delta_{\hbar} e_{i}&=&\sum  e_{k}\otrh e_{l} \cdoth g_{i}^{kl}
+\hbar \tilde{\Delta}_{\hbar}^{1}e_{i}, \label{eq:dhe-tilde}
\end{eqnarray}
for some  $\tilde{\Delta}_{\hbar}^{2}X$ and
$\tilde{\Delta}_{\hbar}^{1}e_{i}\in U_{\hbar}A\otrh
U_{\hbar}A$.

Hence,
\be
(\Delta_{\hbar}\otrh id )\Delta_{\hbar}^{1} X&=&
\sum e_{k} \otrh e_{l}\otrh  (\alpha_{\hbar}  g_{i}^{kl}*_{\hbar}
\beta_{\hbar}
f^{ij}*_{\hbar}e_{j} )+\sum (Y_{i}\otrh 1\otrh Z_{i}
+1\otrh Y_{i}\otrh Z_{i}\\
&&- Z_{i}\otrh 1 \otrh Y_{i}-1\otrh Z_{i}\otrh Y_{i})
+\hbar [\sum  \tilde{\Delta}_{\hbar}^{1} e_{i}\otrh e_{j}
\cdoth f^{ij} \\
&&+\sum \Delta_{\hbar}^{1}Y_{i}\otrh Z_{i}-\sum
\Delta_{\hbar}^{1}Z_{i}\otrh Y_{i}
+(\Delta_{\hbar}\otrh id )\tilde{\Delta}_{\hbar}^{2} X],
\ee

and

\be
\Delta_{\hbar}^{1} X \otrh 1&=&
\sum e_{i} \otrh e_{j}\otrh \alpha_{\hbar}f^{ij}+
\sum_{i}(Y_{i}\otrh Z_{i}\otrh 1-Z_{i}\otrh Y_{i}\otrh 1)+
\hbar \tilde{\Delta}_{\hbar}^{2} X \otrh 1.
\ee

Thus we have  $ J_{\hbar}=I_{\hbar}+\hbar K_{\hbar}$,
where
\be
I_{\hbar} &=& \sum e_{i} \otrh e_{j}\otrh ( \alpha_{\hbar}
g_{l}^{ij}*_{\hbar} \beta_{\hbar} f^{lk}*_{\hbar}e_{k} )
-\sum (\alpha_{\hbar}f^{il}*_{\hbar}\beta_{\hbar} g_{l}^{jk}*_{\hbar}e_{i})
\otrh e_{j}\otrh e_{k} \\
&&-\sum \beta_{\hbar}f^{ij} \otrh e_{i} \otrh e_{j}
+ \sum e_{i} \otrh e_{j}\otrh \alpha_{\hbar}f^{ij},
\ee
and
\be
K_{\hbar}
&=& \sum  \tilde{\Delta}_{\hbar}^{1} e_{i}\otrh e_{j}
\cdoth f^{ij}
+\sum \Delta_{\hbar}^{1}Y_{i}\otrh Z_{i}-\sum \Delta_{\hbar}^{1}Z_{i}\otrh Y_{
i}
+(\Delta_{\hbar}\otrh id )\tilde{\Delta}_{\hbar}^{2} X
+\tilde{\Delta}_{\hbar}^{2} X \otrh 1\\
&&-[\sum f^{ij}\cdoth  e_{i}  \otrh \Delta_{\hbar}^{1} e_{j}
+\sum Y_{i}\otrh \Delta_{\hbar}^{1}Z_{i}-\sum  Z_{i}\otrh \Delta_{\hbar}^{1}Y_
{i}
+( id \otrh  \Delta_{\hbar})\Delta_{\hbar}^{2} X
+1 \otrh \Delta_{\hbar}^{2}X].
\ee

By $\alt $, we denote  the  standard
skew-symmetrization operator on $UA\otr UA\otr UA$:
$$\alt (x_{1}\otr x_{2}\otr x_{3})=\sum_{\sigma\in S_{3}}\frac{1}{3!}
(-1)^{|\sigma |}  x_{\sigma (1)}\otr  x_{\sigma (2)}\otr x_{\sigma (3)}, $$
where $x_{1}, x_{2}, x_{3}\in UA$.
It is tedious but straightforward to
verify that  $\alt  ( \limh \frac{1}{\hbar}I_{\hbar})=0$.
Therefore, $\alt  ( \limh K_{\hbar})=0$, i.e.,
\be
&& \alt [ \sum  f^{ij}\tilde{\Delta}^{1} e_{i}\otr e_{j}
+\sum \Delta^{1}Y_{i}\otr Z_{i}-\sum \Delta^{1}Z_{i}\otr Y_{
i}
+(\Delta\otr id )\tilde{\Delta}^{2} X
+\tilde{\Delta}^{2} X \otr 1 \\
&&-(\sum f^{ij} e_{i}  \otr \Delta^{1} e_{j}
+\sum Y_{i}\otr \Delta^{1}Z_{i}-\sum  Z_{i}\otr \Delta^{1}Y_
{i}
+( id \otr  \Delta)\Delta^{2} X +1 \otr \Delta^{2}X)] =0.
\ee

The  final conclusion  thus follows immediately
by  applying the skew-symmetrization
operator $\alt$ to the equation above and using  the
following simple facts:

\begin{lem}
\begin{enumerate}
\item $\alt (\tilde{\Delta}^{2}X\otr 1-1\otr \Delta^{2}X)=0.$
\item $ \alt  \sum ( f^{ij}\tilde{\Delta}^{1}e_{i}\otr e_{j}-
f^{ij}e_{i}\otr \Delta^{1}e_{j})=0.$
\item $\alt ((id\otr \Delta )\Delta^{2}X)=\alt ( (\Delta \otr id)
\tilde{\Delta}^{2}X)=0.$
\end{enumerate}
\end{lem}
\pf It follows from Equations (\ref{eq:d1h}) and (\ref{eq:dlh-tilde}) that
\be
\hbar {\Delta}^{2}_{\hbar}X- \hbar \tilde{\Delta}^{2}_{\hbar}X&=&
\sum (f^{ij}\cdoth e_{i}\otrh e_{j}-e_{i}\otrh e_{j}\cdoth  f^{ij})\\
&=&\sum( f^{ij}\cdoth e_{i}\otrh e_{j}-e_{i}\otrh f^{ij} \cdoth e_{j}
+e_{i}\otrh f^{ij} \cdoth e_{j} -e_{i}\otrh e_{j}\cdoth  f^{ij})\\
&=&\hbar  \sum (\delta f^{ij} *_{\hbar} e_{i} )\otrh e_{j}
+\hbar  \sum e_{i}\otrh (\delta f^{ij} *_{\hbar} e_{j}) +\oh.
\ee
Hence,
\be
&&{\Delta}^{2}X-\tilde{\Delta}^{2}X\\
&=&\sum (\delta f^{ij})e_{i}\otr e_{j}+\sum e_{i}\otr (\delta f^{ij} )e_{j}\\
&=&\sum \Delta (\delta f^{ij})  e_{i}\otr e_{j},
\ee
which is symmetric.
It thus follows that
\be
&& \alt (\tilde{\Delta}^{2}X\otr 1-1\otr \Delta^{2}X)\\
&=&\alt (\Delta^{2}X\otr 1-1\otr \Delta^{2}X-\sum \Delta
(\delta f^{ij})  e_{i}\otr e_{j}\otr 1)\\
&=&0.
\ee

Similarly, one can show that
$${\Delta}^{1}e_{i}-\tilde{\Delta}^{1}e_{i}=
\sum_{kl} \Delta (\delta g_{i}^{kl})(e_{k }\otr e_{l}).$$
Thus, (ii) follows immediately.
Finally, (iii) is obvious since $\Delta$ is co-commutative. \qed

Combining   Propositions \ref{thm:delta0}-\ref{pro:delta2},
we conclude that $(\oplus \gm (\wedge^{*}A) , \wedge , [\cdot , \cdot ],
\delta
)$
is  indeed a strong differential  Gerstenhaber algebra. Hence $(A, A^* )$ is
a Lie bialgebroid,  which  is
called
the classical limit of the  quantum groupoid  $U_{\hbar}A$.
In summary, we have proved

\begin{thm}
   A  quantum groupoid  $(U_{\hbar}A, R_{\hbar},
   \alpha_{\hbar}, \beta_{\hbar}, m_{\hbar}, \Delta_{\hbar},
\epsilon_{\hbar})$
naturally induces a Lie bialgebroid
$(A, A^* )$ as a classical limit. The   induced Poisson
structure  of  this Lie bialgebroid on the base  manifold  $P$  coincides
with the one obtained as the  classical limit of the base
   $*$-algebra $R_{\hbar}$.
\end{thm}

As an example, in what follows,  we will examine the case where
   the quantum groupoids are obtained from the standard Hopf algebroid
$UA[[\hbar ]]$ by  twists.  Consider $(UA[[\hbar ]],  R[[\hbar ]],
   \alpha, \beta, m, \Delta, \epsilon )$ equipped with the
standard  Hopf algebroid structure induced from that on $UA$.
Assume that
\begin{equation}
\label{eq:twistf}
\calfh =1\otr 1+\hbar \fl + \oh \in UA\otr UA [[\hbar ]],
\end{equation}
where $\fl\in UA\otr UA$,
is a twistor,
   and let $(U_{\hbar}A,  R_{\hbar},
   \alpha_{\hbar}, \beta_{\hbar}, m_{\hbar}, \Delta_{\hbar},
\epsilon_{\hbar})$
be
the resulting twisted QUE algebroid.

\begin{lem}
\label{lem:twist-f}
Assume that $\calfh \in UA\otr UA [[\hbar ]]$ given by
Equation (\ref{eq:twistf}) is a twistor. Then
\begin{equation}
\label{eq:f1}
\fl \otr 1+(\Delta\otr id) \fl
=1 \otr \fl +(id \otr \Delta)  \fl.
\end{equation}
\end{lem}
\pf This follows immediately from  
computing the $\hbar^{1}$-term
in the $\hbar$-expansion
of Equation  (\ref{eq:cocycle}). \qed

Using Proposition \ref{pro:appendix} in the Appendix, we have

\begin{cor}
Under the same hypothesis as in Lemma  \ref{lem:twist-f},  then
$\Lambda =\alt \fl ( \stackrel{def}{=}\fl-\fl_{21})$
is a section of $\wedge^{2}A$.
\end{cor}

Now it is natural to expect the following:

\begin{thm}
\label{thm:twist-limit}
   Let $(U_{\hbar}A,  R_{\hbar},
   \alpha_{\hbar}, \beta_{\hbar}, m_{\hbar}, \Delta_{\hbar},
\epsilon_{\hbar})$
be   the   quantum groupoid obtained from \\
$(UA[[\hbar ]],  R[[\hbar ]],
   \alpha, \beta, m, \Delta, \epsilon )$ by
   twisting via $\calfh$, where
$$\calfh =1\otr 1+\hbar \fl + \oh \in UA\otr UA [[\hbar ]].$$
 Then its classical limit
is a coboundary Lie bialgebroid $(A, A^{*}, \Lambda )$,
where $\Lambda =\alt \fl$.  In particular, its  induced  Poisson   structure
on the base  manifold  is $\rho \Lambda$,  which admits  $R_{\hbar}$  as
a  deformation quantization.
\end{thm}
\pf It suffices to prove that $\delta f=[f, \Lambda ]$
and $\delta  X=[X, \Lambda ],  \forall f\in R $ and $X\in \gm (A)$.

   Write
$$\fl =\sum_{i} d_{i}\otr e_{i}\in UA\otr UA.$$
Using  Equations (\ref{eq:alpha})-(\ref{eq:beta}),
it is easy to see that for any $f\in R$,
\be
\alphah f&=&f +\hbar \sum_{i} ((\rho d_{i})f)e_{i}+\oh , \ \ \mbox{ and }\\
\beta_{\hbar} f&=&f +\hbar \sum_{i} ((\rho e_{i})f)d_{i}+\oh .
\ee
It thus follows that
$$\delta f=\sum_{i}[((\rho d_{i})f)e_{i}-((\rho e_{i})f)d_{i}]
=[f, \Lambda ].$$

Now we will prove the second identity $\delta X=[X, \Lambda ]$. From
the definition of  $\calfh^{\#}$,  it follows that

\be
&&\calfh^{\#}[1\otrh X+X\otrh 1+\hbar  \sum_{i}
(Xd_{i}\otrh e_{i}+d_{i}\otrh Xe_{i}
-d_{i}X\otrh e_{i} -d_{i} \otrh  e_{i}X)]\\
&=&\calfh  [1\ot X+X\ot 1+\hbar \sum_{i} (Xd_{i}\ot e_{i}+d_{i}\ot Xe_{i}
-d_{i}X\ot e_{i} -d_{i} \ot  e_{i}X)]\\
&=&1\otr X+X\otr 1 +\hbar  \sum_{i}(Xd_{i}\otr e_{i}+d_{i}\otr e_{i}X)+\oh\\
&=&(\Delta X)\calfh +\oh.
\ee

It thus follows that
\be
&& \Delta_{\hbar}X\\
&=&{\calfh^{\#}}^{-1} ((\Delta X) \calfh ) \\
&=& 1\otrh X+X\otrh 1+\hbar  \sum_{i} (Xd_{i}\otrh e_{i}+d_{i}\otrh Xe_{i}
-d_{i}X\otrh e_{i} -d_{i} \otrh  e_{i}X)+\oh .
\ee

Therefore
$$\Delta^{1}X= (1\ot X+X\ot 1)\fl -\fl (1\ot X+X\ot 1). $$

This  immediately implies that

\be
\delta X&=&\Delta^{1}X-(\Delta^{1}X)_{21}\\
& =&[(1\ot X+X\ot 1)\fl -\fl (1\ot X+X\ot 1)]-
[(1\ot X+X\ot 1)\fl_{21} -\fl_{21} (1\ot X+X\ot 1)]\\
& =&(1\ot X+X\ot 1)\Lambda -\Lambda (1\ot X+X\ot 1)\\
& =&[X, \Lambda ].
\ee
This concludes the proof. \qed

More generally, we have

\begin{thm}
\label{thm:limit-rmatriod}
Assume that $(U_{\hbar}A,  R_{\hbar}, \alpha_{\hbar}, \beta_{\hbar},
m_{\hbar},
   \Delta_{\hbar}, \epsilon_{\hbar})$
   is a  quantum groupoid with classical limit
$(A, A^* )$. Let $\calfh\in U_{\hbar}A\otrh U_{\hbar}A$ be a twistor such
that $\calfh =1\otrh 1  (mod \hbar )$.
Then $\Lambda = Alt(\limh \hbar^{-1}  (\calfh -1 \otrh 1 ) )$ is
a section of $\wedge^{2} A $,   and is a Hamiltonian
operator of the Lie bialgebroid $(A, A^* )$. If $(U_{\hbar}A,
   \tilde{R}_{\hbar}, \tilde{ \alpha}_{\hbar},  \tilde{\beta}_{\hbar},
m_{\hbar},
   \tilde{\Delta}_{\hbar}, \epsilon_{\hbar})$ is obtained from
$(U_{\hbar}A,  R_{\hbar}, \alpha_{\hbar}, \beta_{\hbar}, m_{\hbar},
   \Delta_{\hbar}, \epsilon_{\hbar})$ by twisting via $\calfh$,
its  corresponding Lie bialgebroid
is obtained from $(A, A^* )$ by twisting via $\Lambda$.
In particular, if $(A, A^* )$ is a coboundary Lie bialgebroid
with $r$-\matroid\   $\Lambda_{0}$, the latter is still
    a coboundary Lie bialgebroid
and   its  $r$-\matroid\  is    $\Lambda_{0}+\Lambda$.
\end{thm}
\pf The proof is similar to that of Theorem \ref{thm:twist-limit},
even though it is a little bit more complicated. We omit it here. \qed
{\bf  Remark.} It is easy to see that  the classical limit of
the quantum groupoid $(\dhh , R_{\hbar},
   \alpha_{\hbar} , \beta_{\hbar} , m, \Delta_{\hbar}, \epsilon )$
in Example 4.1  is the standard  Lie bialgebroid $(TP, T^{*}P)$
associated to  a  Poisson manifold $P$ \cite{MX94}. 
It would be interesting to   explore its ``dual"  quantum
groupoid, namely the one with the Lie bialgebroid $(T^{*}P, TP)$
as its  classical limit.

\section{Quantization of Lie bialgebroids}
\begin{defi}
A quantization of  a  Lie bialgebroid  $(A, A^{*})$ is
a quantum groupoid
   $(U_{\hbar}A,  R_{\hbar}, \alpha_{\hbar}, \beta_{\hbar},\\
   m_{\hbar}, \Delta_{\hbar}, \epsilon_{\hbar})$
   whose classical limit is  $(A, A^{*})$.
\end{defi}

It is a deep theorem of Etingof and Kazhdan \cite{EK} that
every Lie bialgebra is quantizable.
   On the other  hand, the existence of $*$-products for an
arbitrary Poisson manifold was recently proved  by Kontsevich \cite{K}.
   In terms of Hopf algebroids, this amounts  to saying
that the Lie bialgebroid $(TP, T^{*}P)$ associated to
a Poisson manifold $P$  is always quantizable.
It is therefore   natural to expect:
\begin{quote}
{\bf Conjecture}  Every Lie bialgebroid is quantizable.
\end{quote}

Below we will prove a very   special case of this
conjecture by using Fedosov  quantization method  \cite{Fedosov} \cite{Xu0}.

\begin{thm}
\label{thm:quan}
Any regular triangular Lie bialgebroid is quantizable.
\end{thm}

We need some preparation first. Recall that
given  a Lie algebroid  $A\lon P$  with anchor $\rho $,
an $A$-connection on a vector bundle $E\lon P$ is
   an $\reals$-linear  map:
$$\gm (A)\otimes \gm (E)\lon \gm (E)$$
$$X\otimes s\lon \nabla_{X}s, $$
satisfying the axioms resembling   those of  usual
   linear connections, i.e.,
$\forall f\in C^{\infty}(P), \ X\in \gm (A), s\in \gm (E)$,
\be
&&\nabla_{fX}s=f\nabla_{X}s;\\
&&\nabla_{X}(fs)=(\rho (X)f) s+f \nabla_{X}s .
\ee
In particular, if $E=A$, an $A$-connection is called torsion-free if
$$\nabla_{X}Y-\nabla_{Y}X=[X, Y], \ \ \ \forall X, Y\in \gm (A).$$
A torsion-free connection always exists for any Lie algebroid.

Let $\omega\in \gm (\wedge^2 A^* )$ be a closed two-form, i.e.,
$d\omega =0$.
 An $A$-connection on $A$ is said to be compatible
with $\omega$ if $\nabla_{X}\omega =0$,
$\forall X\in \gm (A)$.  If $\omega$ is non-degenerate,
a compatible  torsion-free connection always exists.

\begin{lem}
\label{lem:connection}
If $\omega\in \gm (\wedge^2 A^* )$ is a closed non-degenerate
two-form, there exists a compatible torsion-free $A$-connection on $A$.
\end{lem}
\pf  This result is standard (see \cite{NT} \cite{Weinstein1}).
The proof is simply
a repetition of that of the existence
of a symplectic connection for a symplectic manifold.
For completeness, we sketch a proof here.

First, take any torsion-free $A$-connection $\nabla $. Then
any other $A$-connection can be written as
\begin{equation}
\tilde{\nabla}_{X}Y=\nabla_{X}Y +S(X, Y) , \ \ \ \forall X, Y\in \gm (A),
\end{equation}
where $S$ is a $(2, 1)$-tensor.
Clearly, $\tilde{\nabla}$ is torsion-free if and only if
$S$ is symmetric, i.e., $S(X, Y)=S(Y, X)$ for any $X, Y\in \gm (A)$.

$\tilde{\nabla}$ is compatible with $\omega\in \gm (\wedge^2 A^* )$
   if and only if  $\tilde{\nabla}_{X}\omega =0$.
The latter is equivalent to
\begin{equation}
\omega (S(X, Y), Z)-\omega (S(X, Z), Y)=(\nabla_{X}\omega )(Y, Z).
\end{equation}

Let $S$ be the $(2, 1)$-tensor defined by
the equation:

\begin{equation}
\label{eq:s}
\omega (S(X, Y), Z)=\frac{1}{3}[(\nabla_{X}\omega )(Y, Z)+
(\nabla_{Y}\omega )(X, Z)] .
\end{equation}

Clearly, $S(X, Y)$,   defined  in this way,  is symmetric
with respect to $X$ and $Y$.  Now
\be
&&\omega (S(X, Y), Z)-\omega (S(X, Z), Y)\\
&=&\frac{1}{3}[(\nabla_{X}\omega )(Y, Z)+ (\nabla_{Y}\omega )(X, Z)]
-\frac{1}{3}[(\nabla_{X}\omega )(Z, Y)+(\nabla_{Z}\omega )(X, Y)]\\
&=&\frac{1}{3}[(\nabla_{X}\omega )(Y, Z)+(\nabla_{Y}\omega )(X, Z)
+(\nabla_{X}\omega )( Y, Z)+(\nabla_{Z}\omega )(Y, X)]\\
&=&(\nabla_{X}\omega ) (Y, Z).
\ee
Here the last step follows from  the identity:
$$(\nabla_{X}\omega )(Y, Z)+(  \nabla_{Y}\omega )(Z, X)+
(\nabla_{Z}\omega )(X, Y)=0,$$
which is equivalent to
$d\omega =0$.
This implies that $\tilde{\nabla}$ is a torsion-free symplectic connection.
\qed
{\bf Proof of Theorem \ref{thm:quan}}
   Let $(A, A^*, \Lambda )$ be a regular triangular Lie bialgebroid.
Then  $\Lambda^{\#}: A^* \lon A$ is  a Lie algebroid
morphism \cite{MX94}.
Therefore its image $\Lambda^{\#}A^*$ is a  Lie subalgebroid
of $A$, and $ \Lambda$ can be considered as a section of
$\wedge^{2} (\Lambda^{\#}A^* )$.
Hence, by restricting to $\Lambda^{\#}A^*$ if necessary,
one may always assume that $\Lambda $ is nondegenerate.
Let $\omega =\Lambda^{-1}\in \gm (\wedge^{2} A^* )$. Then
$\omega$ is closed: $d\omega =0$. Let $\nabla$ be  a compatible
torsion-free $A$-connection  on $A$, which always exists
according to Lemma \ref{lem:connection}. Let $\poiddd{G}{P}{}$ be
a local Lie  groupoid corresponding to  the Lie algebroid $A$.
Let $\Lambda^{l}$ denote  the left translation of $\Lambda$,
so $\Lambda^{l}$ defines a left invariant Poisson structure on $G$.
This is a regular Poisson structure, whose  symplectic leaves
are simply $\alpha$-fibers. The $A$-connection  $\nabla$
induces a fiberwise linear connection $\tilde{\nabla}$
   for the $\alpha$-fibrations. To see this, simply
define for any $X, Y\in \gm (A)$,
\begin{equation}
\tilde{\nabla}_{X^{l}}Y^{l}=(\nabla_X Y)^{l},
\end{equation}
   where $X^{l}, Y^{l}$ and $(\nabla_X Y)^{l}$ denote their corresponding
   left invariant vector fields on $G$. Since left invariant
vector fields span the tangent space of $\alpha$-fibers, this
indeed defines a linear connection on each $\alpha$-fiber $\alpha^{-1}(u),
\forall u\in P$,  which  is  
   denoted by $\tilde{\nabla}_u$. Clearly, $\tilde{\nabla}_u$
is torsion-free since $\nabla$ is torsion-free.
Moreover, $\tilde{\nabla}$ preserves the Poisson structure $\Lambda^{l}$,
and  is left-invariant in the sense that
\begin{equation}
L_{x}^{*}\tilde{\nabla}_u=\tilde{\nabla}_v, \  \ \forall
x\in G \ \mbox{ such that } \beta (x)=u , \alpha (x)=v.
\end{equation}

Applying Fedosov quantization  to this situation, one obtains
a $*$-product on $G$:
\begin{equation}
f*_{\hbar} g=fg +\half \hbar \Lambda^{l}(f, g)+\cdots +\hbar^{k}
B_{k}(f, g)+\cdots
\end{equation}
quantizing the  Poisson structure $\Lambda^l$. In fact, this $*$-product is
   given by a  family of leafwise  $*$-products indexed by $u\in P$
quantizing  the  leafwise symplectic structures on $\alpha$-fibers.
The Poisson structure $\Lambda^{l}$ is left invariant, 
so  the  leafwise symplectic structures are invariant under left
translations. Moreover, since  the symplectic connections $\tilde{\nabla}_u$
are left-invariant, the resulting Fedosov $*$-products are invariant under
left translations.
In other words, the  bidifferential operators $B_{k}(\cdot , \cdot )$
   are all left invariant, and therefore
can be considered as elements in  $UA\otr UA $.
 In this way, we  obtain   a formal power series 
$\calf_{\hbar}=1+\half \hbar \Lambda +O(\hbar^2 )\in UA\otr UA[[\hbar ]]$
so  that the $*$-product on $G$ is 
$$f*_{\hbar} g=\calf_{\hbar} (f, g), \ \ \ \ \forall f, g\in C^{\infty}(G).$$
The associativity of $*_{\hbar}$ implies that $\calf_{\hbar}$ satisfies
Equation (\ref{eq:cocycle}):
\begin{equation}
   (\Delta \otr  id )\calf_{\hbar}   \calf_{\hbar}^{12}
~ = ~ (id \otr  \Delta) \calf_{\hbar}  \calf_{\hbar}^{23}.
\end{equation}
The identity $1*_{\hbar}f=f*_{\hbar}1=f$ implies that
\begin{equation}
(\epsilon \otr id) \calf_{\hbar} ~ = ~ 1_{H}; \ \
(id \otr \epsilon) \calf_{\hbar} ~ = ~ 1_{H}.
\end{equation}
Thus $\calf_{\hbar}\in UA\otr UA [[\hbar ]]$ is  a twistor, and
the resulting  twisted Hopf algebroid
$(U_{\hbar}A, R_{\hbar}, \alpha_{\hbar}, \beta_{\hbar}, m_{\hbar},
   \Delta_{\hbar}, \epsilon_{\hbar})$
is a quantization of the triangular Lie bialgebroid  $(A, A^*, \Lambda )$
according to Theorem \ref{thm:twist-limit}.
This concludes the proof of the theorem. \qed

In particular, when the base $P$ reduces to a point,  Theorem \ref{thm:quan}
implies  that every finite dimensional triangular
$r$-matrix is quantizable. Of course, there is no need to use Fedosov
method in  this   case.  There is  a very nice
short proof due to Drinfel'd \cite{Drinfeld2}.\\\\

We note that the induced  Poisson structure
on  the base   manifold of a non-degenerate  triangular Lie bialgebroid,
also called {\em a symplectic Lie algebroid},
was  studied by Nest-Tsygan \cite{NT} and
Weinstein \cite{Weinstein1}, for which a $*$-product was constructed.
   Indeed, our algebra $R_{\hbar}$  provides
a  $*$-product for such a  Poisson structure, where
the multiplication is simply defined by the
push forward of  $\calf_{\hbar}$ under the anchor $\rho$:
$$a*_{\hbar}b=(\rho \calf_{\hbar})(a, b), \ \ \ \forall a, b
\in C^{\infty}(P)[[\hbar ]].$$
 So here we obtain an
alternative proof of (a slightly more general version of) 
their quantization result.

\begin{cor}
The induced  Poisson structure
on  the base   manifold of a regular (in particular non-degenerate)  triangular Lie bialgebroid
is quantizable.
\end{cor}

\section{Dynamical quantum groupoids}

This section is devoted to the study of an important example of
quantum groupoids, which are connected with the so called
quantum dynamical $R$-matrices.
Let  $\uqg$ be   a quasi-triangular quantum  universal enveloping
algebra over $\complex$ with $R$-matrix $R\in \uqg \ot \uqg$,
$\eta \subset \frakg$ a finite  dimensional Abelian Lie  subalgebra
such that $U\eta [[\hbar ]]$ is  a commutative subalgebra of $\uqg$.
By  $\C (\eta^{*})$, we denote  the algebra of meromorphic functions
   on $\eta^{*}$, and  by $\cald$ the algebra
of meromorphic differential operators on $\eta^{*}$.
Consider $H=\cald\ot \uqg $.
Then $H$ is a Hopf algebroid  over $\complex$ with base algebra
$R=\C (\eta^{*})[[\hbar ]]$,   whose
coproduct and counit  are denoted, respectively, by
$\Delta$ and $\epsilon$. Moreover the map
\begin{equation}
\mu (D\otimes u)(f)=(\epsilon_{0} u) D(f), \ \ \forall D\in \cald ,
\ u\in \uqg , \  f\in \C(\eta^{*})[[\hbar ]],
\end{equation}
is an anchor map. Here $\epsilon_0$ is the counit  of the Hopf algebra
   $\uqg$. Let us   fix a basis in  $\eta$, say $\{h_{1},  \cdots ,h_{k}\}$,
and let $\{\xi_{1}, \cdots , \xi_{k}\}$ be its dual basis,
which in turn
    defines a coordinate system $(\lambda_{1}, \cdots ,\lambda_{k})$
on $\eta^{*}$.

Set
\begin{equation}
\label{eq:theta}
\theta =\sum_{i=1}^{k}(\frac{\partial}{\partial \lambda_{i}}\ot h_{i})
\in H\ot H , \ \mbox{ and } \Theta =\exp{\hbar \theta} \in H\ot H.
\end{equation}
Note that $\theta $, and hence $\Theta$, is independent   of
the choice of  bases in  $\eta$.
The following fact  can be easily verified.

\begin{lem}
\label{lem:theta}
$\Theta $ satisfies Equations (\ref{eq:cocycle}) and
(\ref{eq:co-unit}).
\end{lem}
\pf Consider
   $H_{0}=\cald^{inv}\ot U\eta [[\hbar ]]$, where $\cald^{inv}$
consists of holomorphic  differential operators on $\eta^{*}$ invariant under
the translations.
Then $H_{0}$ is a Hopf subalgebroid of $H$, which is in fact a Hopf algebra.
Clearly, $\theta \in H_{0}\ot H_{0}$, so $\Theta \in H_{0}\ot H_{0}$. It thus
suffices to prove that

$$(\Delta \ot  id )\Theta  \Theta^{12}
=(id \ot  \Delta) \Theta \Theta^{23}  \ \ \mbox{ in } H_{0}\ot H_{0}\ot H_{0}.$$

Now
\be
&&(\Delta \ot  id )\Theta  \Theta^{12}\\
&=&((\Delta \ot  id ) \exp{\hbar \theta} )  \exp{\hbar \theta^{12}}\\
&=&\exp\hbar ((\Delta \ot  id )\theta + \theta^{12})\\
&=&\exp\hbar \sum_{i=1}^{k}( \parr{}\ot 1\ot h_{i}+1\ot \parr{}\ot h_{i}
+\parr{}\ot h_{i}\ot 1 	).
\ee
Here in  the second equality we
    used the fact that $(\Delta \ot  id ) \theta$ and
$\theta^{12}$ commute in  $H_{0}\ot H_{0}\ot H_{0}$.

Similarly, we  have
$$(id \ot  \Delta) \Theta \Theta^{23}=\exp\hbar \sum_{i=1}^{k}( \parr{}\ot
1\ot h_{i}+1\ot \parr{}\ot h_{i}
+\parr{}\ot h_{i}\ot 1).$$

This proves Equation (\ref{eq:cocycle}).
Finally, Equation (\ref{eq:co-unit})
   follows from a straightforward verification. \qed
{\bf  Remark.} There is a more intrinsic way of proving this
fact in terms of deformation quantization.
Consider  $T^{*}\eta^{*}$ equipped with the standard cotangent
bundle symplectic structure $\sum_{i=1}^{k} d\lambda_{i}\wedge dp_{i}$.
It is  well-known that,  $\forall f, g\in C^{\infty}(T^{*}\eta^{*})[[\hbar ]]$,
\begin{equation}
f*_{\hbar}g=fe^{\hbar (\sum_{i=1}^{k} \overleftarrow{\parr{}}\ot
\overrightarrow{\parrr{}{p_{i}}})}g
\end{equation}
defines a $*$-product on $T^{*}\eta^{*}$,
   called the Wick type  $*$-product  corresponding to the
normal ordering  quantization.  Hence
$\exp {\hbar (\sum_{i=1}^{k} \overleftarrow{\parr{}}\ot
\overrightarrow{\parrr{}{p_{i}}})}$,
as  a (formal) bidifferential operator on $T^{*}\eta^{*}$
(i.e. as  an element in $\cald (T^{*}\eta^{*} )\ot_{\C (T^{*}\eta^{*})}
\cald (T^{*}\eta^{*}  )[[\hbar ]] $)
satisfies Equations (\ref{eq:cocycle})-(\ref{eq:co-unit})
according to Example 4.1.
   Note that elements
in $(\cald  \ot U\eta ) \ot_{\C(\eta^{*})}(\cald  \ot U\eta )
[[\hbar ]]$
can be considered as  (formal) bidifferential operators on $T^{*}\eta^{*}$
invariant under the  $p$-translations, so
$(\cald  \ot U\eta ) \ot_{\C(\eta^{*})}(\cald  \ot U\eta
)[[\hbar ]]$
is naturally  a subspace of
   $\cald (T^{*}\eta^{*} )\ot_{\C(T^{*}\eta^{*})} \cald (T^{*}\eta^{*}
)[[\hbar ]]$.
Clearly $\exp {\hbar (\sum_{i=1}^{k} \overleftarrow{\parr{}}\ot
   \overrightarrow{ \parrr{}{p_{i}}})}$
is   a $p$-invariant bidifferential operator
on $T^{*}\eta^{*}$, and is equal
   to $\Theta $ under the above identification.
Equations (\ref{eq:cocycle}) and (\ref{eq:co-unit}) thus
follow immediately.\\\\


In other words, $\Theta$ is  a twistor
   of  the Hopf algebroid $H$. As we  see below,
it is this $\Theta$ that links  a
shifted cocycle $F(\lambda)$ and   a Hopf algebroid twistor.

Given  $F\in \C(\eta^{*}, \uqg\ot \uqg)$, define
$F^{12} (\lambda +\hbar h^{(3)})\in \C(\eta^{*},  \uqg\ot \uqg \ot \uqg )$
by
\begin{eqnarray}
F^{12} (\lambda +\hbar h^{(3)}) &=&F(\lambda  )\ot 1+\hbar \sum_{i}
\parr{F}\ot h_{i}
+\frac{1}{2!}\hbar^{2}\sum_{i_{1}i_{2}}
\frac{\partial^{2}F}{\partial \lambda_{i_{1}}\partial \lambda_{i_{2}}}\ot
h_{i_{1}}h_{i_{2}} \nonumber\\
&& \ \ \ \ +\cdots
+\frac{\hbar^k}{k!}\sum \frac{\partial^{k} F}{\partial  \lambda_{i_{1}}
\cdots \partial \lambda_{i_{k}}} \ot h_{i_{1}}\cdots h_{i_{k}} +\cdots,
\end{eqnarray}
similarly for $F^{23} (\lambda +\hbar h^{(1)})$ etc.

\begin{lem}
Let  $X$ be a
meromorphic vector field on $\eta^*$,
and $F\in \C(\eta^{*}, \uqg\ot \uqg )$. Then
\begin{equation}
\label{eq:delxf}
(\Delta X) F=F(\Delta X) +X(F), \ \ \ \ \mbox{in }\  H\otr H.
\end{equation}
\end{lem}
\pf Note that $F$, being considered as an element
in $H\otr H$, clearly
satisfies the condition that $\forall f\in R$,
$$F(f\ot 1-1\ot f)=0, \ \ \ \mbox{ in }  H\otr H.  $$
So both $(\Delta  X)  F$ and $F(\Delta X)$
are well defined elements  in $H\otr H$.
For any $f, g \in R$,  considering  both
sides of Equation (\ref{eq:delxf}) as $\uqg\ot \uqg$-valued bidifferential
operators and  applying  them to $f\ot g$, one obtains that
\be
&&[ (\Delta X)F -F(\Delta X)](f\ot g)\\
&=&X(Ffg)-FX(fg)\\
&=&X(F)fg.
\ee
Thus Equation (\ref{eq:delxf}) follows. \qed

An element $F\in \C(\eta^{*}, \uqg \ot \uqg)$ is said to be of zero weight  if

\begin{equation}
\label{eq:0weight}
[F(\lambda ), 1\ot h+h\ot 1]=0, \ \forall \lambda \in \eta^{*}, \ h\in \eta.
\end{equation}

\begin{lem}
\label{lem:shift}
Assume that $F\in \C(\eta^{*}, \uqg \ot \uqg)$ is of zero weight. Then
$\forall
n\in \Bbb{N}$,
\begin{eqnarray}
&&[(\Delta \otr  id )\theta^{n}]F^{12}
=\sum_{k=0}^{n}\sum_{0\leq i_{1} ,\cdots , i_{k} \leq n}
{C_{n}^{k}}(\frac{\partial^{k}F}{\partial \lambda_{i_{1}}\cdots \partial
\lambda_{i_{k}}} \ot h_{i_{1}}\cdots h_{i_{k}} )
(\Delta \otr  id )\theta^{n-k}; \label{eq:ditf}\\
&&[( id \otr \Delta )\theta ] F^{23}=F^{23} ( id \otr \Delta )\theta .
\label{eq:idtf}
\end{eqnarray}
\end{lem}
\pf  To prove Equation (\ref{eq:ditf}),
 let us first
    consider $n=1$.   Then
\be
&&[(\Delta \otr  id )\theta  ]F^{12}\\
&=&\sum_{i} (\Delta \parr{}\otr h_{i}) F^{12}\\
&=&\sum_{i} (\Delta \parr{})F\otr h_{i} \\
&=& \sum_{i} (F\Delta \parr{}+\parrr{F}{\lambda_{i}})\otr h_{i} \\
&=&F^{12}(\Delta \otr  id )\theta +\sum_{i}\parrr{F}{\lambda_{i}}\otr h_{i}.
\ee
The general case follows from   induction,
using the above equation repeatedly.

For Equation (\ref{eq:idtf}), we have
\be
&&[( id \otr \Delta )\theta ] F^{23}\\
&=&\sum_{i}(\parr{}\otr \Delta h_{i})F^{23}\\
&=&\sum_{i}\parr{}\otr (\Delta h_{i} )  F\\
&=&\sum_{i}\parr{}\otr (F \Delta h_{i} )\\
&=&F^{23} ( id \otr \Delta )\theta ,
\ee
where the second  from the last equality follows from
the fact that $F$ is of zero weight. Equation (\ref{eq:idtf}) thus
follows.
\qed

\begin{pro}
\label{pro:shift}
Assume that $F\in \C(\eta^{*}, \uqg \ot \uqg)$ is of zero weight. Then
\begin{eqnarray}
&&[(\Delta \otr  id )\Theta  ]F^{12}(\lambda)
= F^{12} (\lambda +\hbar h^{(3)})(\Delta \otr  id )\Theta;
\label{eq:f_shift}\\
&&[( id \otr \Delta )\Theta ] F^{23}(\lambda)=F^{23}(\lambda)
   ( id \otr \Delta )\Theta  .\label{eq:f-shift}
\end{eqnarray}
\end{pro}
\pf Note that $ (\Delta \otr  id )\Theta
=\exp\hbar ((\Delta \otr  id )\theta )$.
Equations (\ref{eq:f_shift}) and (\ref{eq:f-shift}) thus
   follow  immediately from Lemma  \ref{lem:shift}.  \qed
{\bf Remark. }  One may rewrite Equation (\ref{eq:f_shift}) as
\be
F^{12} (\lambda +\hbar h^{(3)})&=&[(\Delta \otr  id )\Theta ]F^{12}(\lambda )
[(\Delta \otr  id )\Theta ]^{-1}\\
&=&e^{\hbar \sum_{i=1}^{k}(\Delta \parr{}\ot  h_{i})}F^{12}(\lambda )
e^{-\hbar \sum_{i=1}^{k}(\Delta \parr{}\ot  h_{i})}.
\ee
This is essentially the definition of $F^{12} (\lambda +\hbar h^{(3)})$
used in \cite{ABE},
   where the operator $\sum_{i=1}^{k}(\Delta \parr{}\ot  h_{i})$  was
denoted by $\sum_{i=1}^{k}\parr{}h_{i}^{(3)}$.

Now set
\begin{equation}
\F=F(\lambda )\Theta\in H\otr H.
\end{equation}

\begin{thm}
Assume that $F\in \C(\eta^{*}, \uqg \ot \uqg)$ is of zero weight.
Then $\F$ is a twistor (i.e. satisfies Equations (\ref{eq:cocycle})-(\ref{eq:co-unit}))
if and only if
\begin{eqnarray}
\label{eq:shifted}
&&[(\Delta_0  \ot  id )F(\lambda ) ]  F^{12} (\lambda +\hbar h^{(3)})
   = [(id \ot  \Delta_{0} )  F (\lambda )] F^{23}(\lambda  ),\\
&&\label{eq:co}
   (\epsilon_{0} \ot id) F(\lambda )  =  1; \ \
(id \ot \epsilon_{0}) F (\lambda ) =  1,
\end{eqnarray}
where $\Delta_0$ is the coproduct of $\uqg$, and
$\epsilon_{0}$ is the counit map.
\end{thm}
\pf   Using Proposition \ref{pro:shift}, we have
\be
&& \dif\\
&=&[(\Delta \otr  id)F(\lambda )\Theta ]F^{12}(\lambda )\Theta^{12}\\
&=&[(\Delta_{0} \ot id)F(\lambda )][ (\Delta \otr  id) \Theta]
   F^{12}(\lambda )\Theta^{12}\\
&=&[(\Delta_{0} \ot id)F(\lambda )] F^{12} (\lambda +\hbar h^{(3)})
(\Delta \otr  id) \Theta\Theta^{12},
\ee
and
\be
&& \idf\\
&=& [(id \otr  \Delta ) F(\lambda ) \Theta ] F^{23}(\lambda )\Theta^{23}\\
&=&[(id \ot \Delta_{0})F(\lambda ) ](id \otr  \Delta )
   \Theta F^{23}(\lambda )\Theta^{23}\\
&=&[(id \ot \Delta_{0})F(\lambda )]  F^{23}(\lambda )
(id \otr  \Delta )  \Theta \Theta^{23}.
\ee
Thus it follows from Lemma \ref{lem:theta}
   that Equation (\ref{eq:cocycle}) and Equation (\ref{eq:shifted})
are   equivalent.

For Equation (\ref{eq:co-unit}), we note
that $\F =F(\lambda )\Theta =\sum_{k=0}^{\infty}
\frac{\hbar^k}{k!}F(\lambda )\theta^{k}$.
It is easy to see  that  for $k\geq 1$,
    $(\epsilon\otr id)( F(\lambda )\theta^{k})=
(id \otr \epsilon ) (F(\lambda )\theta^{k})=0$
since $\epsilon  (\frac{\partial^{k}}{\partial \lambda_{i_{1}}\cdots
\partial \lambda_{i_{k}}} )=0$
and
$\epsilon (h_{i_{1}}\cdots h_{i_{k}})=0$.
Thus it is immediate that Equations (\ref{eq:co-unit}) and  (\ref{eq:co})
   are equivalent.
This concludes the proof of the theorem.
\qed

A solution to Equations (\ref{eq:shifted})-(\ref{eq:co}) is often called {\em a
shifted cocycle}  \cite{Arnaudon} \cite{BBB} \cite{Jimbo}.
 Moreover, if  $\uqg$ is a quasi-triangular  Hopf algebra
with a quantum   $R$-matrix $R$ satisfying the quantum Yang-Baxter
equation, then 
$R(\lambda )=F^{21}(\lambda )^{-1}R F^{12}(\lambda )$ is a 
solution of the {\em  quantum dynamical Yang-Baxter equation} \cite{BBB}:
\begin{equation}
R^{12}(\lambda )R^{13}(\lambda +\hbar h^{(2)} )R^{23}(\lambda )
=R^{23}(\lambda +\hbar h^{(1)} )R^{13}(\lambda )R^{12}(\lambda +\hbar h^{(3)} )
.
\end{equation}

Now assume that $F(\lambda )$ is  a solution to Equations
   (\ref{eq:shifted})-(\ref{eq:co})
so that we can form a quantum groupoid by twisting
$\cald\ot \uqg $ via $\F$.
The resulting
   quantum groupoid is denoted by $\cald\oth \uqg$,
and is called a {\em dynamical quantum groupoid}.

As an immediate consequence of  Theorem \ref{thm:twist},
   we have

\begin{thm}
As a monoidal category, the category of   $\cald\oth\uqg$-modules
is equivalent to  that of $\cald \ot \uqg $-modules,
and therefore is a braided monoidal category.
\end{thm}
{\bf Remark.}
It is expected that    representations of a quantum dynamical $R$-matrix \cite{Felder}
 can be understood
using this   monoidal category of $\cald \oth \uqg $-modules.
The further
relations between these two objects  will be investigated elsewhere.

In what follows, we   describe various structures
of   $\cald\ot_{\hbar }\uqg$ more  explicitly.

\begin{pro}
\begin{enumerate}
\item $f*_{\F}g=fg, \ \ \forall f, g\in \C(\eta^{*})[[\hbar ]]$, i.e., $R_{\F}$ is
the usual algebra of  functions.
\item  $\alpha_{\F}f= \exp{(\hbar
\sum_{i=1}^{k}h_{i}\frac{\partial}{\partial \lambda_{i}}})f=\sum_{1\leq
i_{1}, \cdots  ,i_{n}\leq k}
\frac{\hbar^{n}}{n!}
   \frac{\partial^{n}f}{\partial \lambda_{i_{1}}\cdots
\partial\lambda_{i_{n}}}h_{i_{1}}\cdots h_{i_{n}}$, $\forall  f \in 
\C(\eta^{*})[[\hbar ]]$;
\item  $\beta_{\F}f=f$,  $\forall  f \in \C(\eta^{*})[[\hbar ]]$.
\end{enumerate}
\end{pro}
\pf Assume that $F(\lambda )=\sum_{i}F_{i}(\lambda )u_{i}\ot v_{i}$,
with $u_{i}, v_{i}\in \uqg$. Let
$$\F_{n}=F (\lambda )\theta^{n}=
\sum_{i} \sum_{1\leq
i_{1}, \cdots  ,i_{n}\leq k}F_{i}(\lambda )u_{i}
\frac{\partial^{n}}{\partial \lambda_{i_{1}}\cdots
\partial\lambda_{i_{n}}} \otr v_{i}h_{i_{1}}\cdots h_{i_{n}}. $$
Then
\be
\alpha_{\F_{n}}f&=&
\sum_{i} \sum_{1\leq i_{1}, \cdots  ,i_{n}\leq k}F_{i}(\lambda )
(\epsilon_{0}u_{i})
\frac{\partial^{n}f}{\partial \lambda_{i_{1}}\cdots
\partial\lambda_{i_{n}}} v_{i} h_{i_{1}}\cdots h_{i_{n}} \\
&=&\sum_{1\leq i_{1}, \cdots  ,i_{n}\leq k}
\frac{\partial^{n}f}{\partial \lambda_{i_{1}}\cdots  \partial
\lambda_{i_{n}}}
[\sum_{i} F_{i}(\lambda ) (\epsilon_{0}u_{i})v_{i}]h_{i_{1}}\cdots h_{i_{n}}\\
&=&\sum_{1\leq i_{1}, \cdots  ,i_{n}\leq k} \frac{\partial^{n}f}{\partial
   \lambda_{i_{1}}\cdots \partial \lambda_{i_{n}}}
h_{i_{1}}\cdots h_{i_{n}},
\ee
where the last equality used the fact that
$\sum_{i} F_{i}(\lambda ) (\epsilon_{0}u_{i})v_{i}=(\epsilon_{0}\ot
id)F(\lambda )=1$.

Similarly, we have $\beta_{\F_{n}}f=f$ if $n=0$, and otherwise
$\beta_{\F_{n}}f
=0$.

Combining these equations, one immediately obtains  that
\be
&&\alpha_{\F}f=\sum_{n=0}^{\infty}
\frac{\hbar^{n}}{n!} \alpha_{\F_{n}}f=
\sum_{1\leq i_{1}, \cdots  ,i_{n}\leq k} \frac{\hbar^{n}}{n!}
   \frac{\partial^{n}f}{\partial \lambda_{i_{1}}\cdots
\partial\lambda_{i_{n}}}h_{i_{1}}\cdots h_{i_{n}},  \\
&&\beta_{\F}f=\sum_{n=0}^{\infty}\frac{\hbar^{n}}{n!} \beta_{\F_{n}}f=f, \ \ \
\mbox{ and }\\
&& f*_{\F}g=(\alpha_{\F}f)(g)=fg.
\ee
\qed

As it is standard \cite{Arnaudon} \cite{Jimbo}, using
$F(\lambda )$,  one  may form a family of  quasi-Hopf algebras
$(\uqg , \Delta_{\lambda})$, where the coproducts
are given by $\Delta_{\lambda}=F(\lambda )^{-1}\Delta_{0}F(\lambda )$.
To describe the relation between  $\Delta_{\F}$
and these quasi-Hopf  coproducts $\Delta_{\lambda}$,
we need to introduce a    ``projection" map from $H\otrf H$ to
   $\C(\eta^{*},  \uqg\ot \uqg
)$.
This can be defined as follows.
Let $Ad_{\Theta}: H\ot H\lon H\ot H$ be the adjoint operator:
$Ad_{\Theta}w=\Theta w \Theta^{-1}$,  $\forall w\in H\ot H$. Composing
with the natural projection, one obtains a
map, denoted by the same symbol $Ad_{\Theta}$,
from $H\ot H$ to $H\otr H$. Since $\alpha_{\Theta}=\alpha_{\F}$,
$\beta_{\Theta}= \beta_{\F}$ and
$\Theta (\beta_{\Theta}f \ot 1-1\ot \alpha_{\Theta} f)=0$,
$\forall f\in R$, in $H\otr H$,
then $\Theta (\beta_{\F}f \ot 1-1\ot \alpha_{\F}f)=0$ in $H\otr H$.
This implies that $Ad_{\Theta}$ descends to a map
from $H\otrf H$ to $H\otr H$. On the other hand,
there exists  an obvious  projection map   $Pr$ from $H\otr H$ to
$\C(\eta^{*}, \uqg\ot \uqg )$, which is just taking
   the 0th-order component.
Now composing  with this projection, one
obtains a map from $H\otrf H$ to $\C(\eta^{*},  \uqg\ot \uqg )$,
which is denoted by $T$.
The following proposition  gives an explicit  description of this map $T$.

An element  $x=D\oth u\in H$, where $D\in \cald [[\hbar ]]$ and
$u\in \uqg$,  is said to be of order $k$ if
   $D$ is a homogeneous differential operator of order $k$.

\begin{pro}
\begin{enumerate}
\item $T(x\otrf y)=0$ if either $x$ or $y$ is
of order greater than zero.
\item $T(fu\otrf gv)
=\sum_{1\leq i_{1}, \cdots  ,i_{n}\leq k} \frac{\hbar^n}{n!}
\frac{\partial^{n}f}{\partial \lambda_{i_{1}}\cdots
\partial\lambda_{i_{n}}}g ( u\ot   h_{i_{1}}\cdots h_{i_{n}}v )
= gu\otr (\alpha_{\F}f) v$, 
$\forall f, g\in \C(\eta^{*})$ and $u, v\in \uqg$.
\end{enumerate}
\end{pro}
\pf  (i) is obvious. We  prove (ii) below.
\be
&&T(fu\otrf gv)\\
&=&Pr (e^{\hbar \sum_{i=1}^{k}(\frac{\partial}{\partial \lambda_{i}}\ot
h_{i})}
(fu\ot gv) e^{-\hbar \sum_{i=1}^{k}(\frac{\partial}{\partial \lambda_{i}}
\ot h_{i})})\\
&=& Pr (e^{\hbar \sum_{i=1}^{k}(\frac{\partial}{\partial \lambda_{i}}\ot
h_{i})}
(fu\ot gv) )\\
&=& Pr  \sum_{1\leq i_{1}, \cdots  ,i_{n}\leq k} \frac{\hbar^n}{n!}
(\frac{\partial^{n} fu}{\partial \lambda_{i_{1}}\cdots
\partial\lambda_{i_{n}}} )\ot g  h_{i_{1}}\cdots h_{i_{n}}v\\
&=&\sum_{1\leq i_{1}, \cdots  ,i_{n}\leq k} \frac{\hbar^n}{n!}
\frac{\partial^{n}f}{\partial \lambda_{i_{1}}\cdots
\partial\lambda_{i_{n}}} g( u\ot   h_{i_{1}}\cdots h_{i_{n}}v)\\
&=&gu\otr (\alpha_{\F}f) v.
\ee
\qed

\begin{pro} The following diagram:
\begin{equation}                         \label{eq:addition}
\matrix{&& \Delta_{\F} &&\cr
          &H  &\vlra&  H\otrf  H&\cr
                  &&&&\cr
                       i&\Bigg\uparrow&&\Bigg\downarrow&T\cr
                               &&&&\cr
                                       &  \uqg  &\vlra&   \C(\eta^{*},
\uqg\ot \uqg )&\cr
                                               &&  \Delta_{\lambda}&&\cr}
\end{equation}
commutes, where $i: \uqg\lon H$ is the natural embedding. I.e.,
$\Delta_{\lambda}=T\smalcirc \Delta_{\F}\smalcirc i$.
\end{pro}
\pf For any $u\in \uqg$, $(\Delta_{\F}\smalcirc i)(u)=\Delta_{\F} (u)
=\F^{-1}(\Delta_{0}u)\F=\Theta^{-1}F(\lambda )^{-1}(\Delta_{0}u )F(\lambda)
\Theta$. Then
$(T\smalcirc \Delta_{\F}\smalcirc i)(u)=T[(\Delta_{\F}\smalcirc i)(u)]=
Pr[F(\lambda )^{-1}(\Delta_{0}u )F(\lambda)]=F(\lambda )^{-1}(\Delta_{0}u
)F(\lambda)
=\Delta_{\lambda} u$.  The conclusion thus follows. \qed

The following theorem describes the classical limit of  the dynamical
quantum groupoid $\cald \oth \uqg$.

\begin{thm}
\label{thm:limit-dug}
Let $(\uqg , R)$ be a quasitriangular quantum universal enveloping
algebra, and $R=1+\hbar r_{0} ( mod \hbar )$.
Assume that $F(\lambda ) \in \uqg\ot \uqg$ is a shifted cocycle
and
that $F(\lambda )=1+\hbar f(\lambda ) ( mod \hbar )$. Then the classical
limit of the corresponding
   dynamical quantum groupoid $\cald \oth \uqg$ is a coboundary
Lie bialgebroid $(A, A^*, \Lambda )$, where $A=T\eta^{*}
   \times \frakg$ and
$\Lambda = \sum_{i=1}^{k}\frac{\partial}{\partial \lambda_{i}}\wedge  h_{i}
+ \alt ( \half r_{0}+  f(\lambda ) )$.
\end{thm}
\pf It is well known \cite{Drinfeld1} that $r_{0}\in \frakg \ot \frakg $,
   and  the operator $\delta: \frakg \lon \wedge^{2}\frakg$,
$\delta a=[1\ot a +a\ot 1 , r_{0}]$, $\forall a\in \frakg $,
defines the cobracket of the corresponding Lie bialgebra of
$(\uqg , R)$. Thus the Lie
bialgebroid corresponding to
$\cald \ot \uqg$ is a coboundary Lie
bialgebroid
   $(T\eta^*  \times \frakg, T^*\eta^*  \times \frakg^* )$
   with $r$-\matroid\  $\half \alt (r_{0})$.
On the other hand, it is obvious
that  $\alt~ \limh \hbar^{-1} (\calf -1) =\alt
f(\lambda )+
\sum_{i=1}^{k}\frac{\partial}{\partial \lambda_{i}}\wedge  h_{i}$.
   The conclusion thus  follows from Theorem \ref{thm:limit-rmatriod}. \qed

As a consequence, we have

\begin{cor}
Under the same
hypotheses as in Theorem \ref{thm:limit-dug},
$r(\lambda )= \alt ( \half r_{0}+ f(\lambda ) )$ is a classical dynamical
$r$-matrix.
\end{cor}

We refer the interested reader to \cite{Arnaudon}  \cite{ESS} \cite{F} 
    \cite{Jimbo}  for
an explicit  construction of  shifted cocycles $F(\lambda )$ for
semisimple Lie algebras.

We end this section by the following:\\\\
{\bf Remark. }  We  may replace
$\theta$  in  Equation (\ref{eq:theta}) by
$\tilde{\theta} =\sum_{i=1}^{k}\half (\frac{\partial}{\partial \lambda_{i}}
\ot h_{i} -h_{i}\ot \frac{\partial}{\partial \lambda_{i}} )
\in H\ot H $ and  set $\tilde{\Theta} =\exp{\hbar \tilde{\theta}} \in H\ot
H$.
It is easy to show that $\tilde{\F}=F(\lambda )\tilde{\Theta }$
satisfies Equation (\ref{eq:cocycle}) is equivalent to
   the following  condition for $F(\lambda )$:
$$[(\Delta_0  \ot  id )F(\lambda )]  F^{12} (\lambda +\half \hbar h^{(3)})
   =  [(id \ot  \Delta_{0} )  F (\lambda )]F^{23}(\lambda  -\half \hbar
h^{(1)}
).$$
In this case, $R(\lambda )=F^{21}(\lambda )^{-1}R F^{12}(\lambda )$ satisfies
the symmetrized quantum  dynamical Yang-Baxter equation:
\begin{equation}
\label{eq:dybe1}
R^{12}(\lambda-\half \hbar h^{(3)})R^{13}(\lambda +\half \hbar h^{(2)} )
R^{23}(\lambda-\half \hbar h^{(1)})
=R^{23}(\lambda  +\half \hbar h^{(1)})R^{13}(\lambda -\half\hbar h^{(2)})
R^{12}(\lambda  + \half \hbar h^{(3)}).
\end{equation}
In fact, both  $\Theta$ and $\tilde{\Theta}$ can be  obtained
from the quantization of the  cotangent bundle
symplectic structure $T^{*}\eta^{*}$, using
the normal
ordering and the Weyl
ordering
respectively, so they are equivalent. This indicates that
solutions to Equation (\ref{eq:dybe}) and Equation (\ref{eq:dybe1})
are   equivalent as well.


\section{Appendix and open questions}

Given any element   $(i_{1}i_{2}i_{3})$ in the symmetric group
$S_{3}$,
by $\sigma_{i_{1}i_{2}i_{3}}$ we denote the  permutation operator on
   $UA\otr UA \otr UA$ given by
$$\sigma_{i_{i}i_{2}i_{3}} ( x_{1}\otr x_{2}\otr x_{3})=
   x_{i_{1}}\otr  x_{i_{2}}\otr x_{i_{3}}. $$

\begin{pro}
\label{pro:appendix}
Assume that $T\in UA\otr UA$ satisfies
\begin{equation}
\label{eq:delta0}
T \otr 1+(\Delta\otr id)T
=1 \otr  T +(id \otr \Delta)  T.
\end{equation}
Then $\alt T\stackrel{def}{=}T-T_{21}$ is a section
of $\wedge^{2}A$.
\end{pro}
\pf First we show that
\label{lem:delta1}
\begin{equation}
\label{eq:delta}
(\Delta \otr id)\alt T=1\otr \alt T+\sigma_{132}(\alt T \otr 1).
\end{equation}

To prove this,
write $T=\sum_{i} u_{i}\otr v_{i}$, where
$u_{i}, v_{i}\in UA$. Then Equation (\ref{eq:delta0}) becomes
\begin{equation}
\label{eq:16}
\sum_{i}u_{i}\otr v_{i}\otr 1+\sum_{i} \Delta u_{i}\otr v_{i}
=\sum_{i} 1\otr u_{i}\otr v_{i}+\sum_{i} u_{i}\otr \Delta v_{i}.
\end{equation}

Applying the permutation operators $\sigma_{231}$ and $\sigma_{132}$,
respectively,  on  both sides of the above equation,
one leads to
\begin{equation}
\label{eq:17}
\sum_{i}v_{i}\otr 1\otr  u_{i}+\sum_{i} \sigma_{231}( \Delta u_{i}\otr v_{i} )
=\sum_{i} u_{i}\otr v_{i} \otr 1+\sum_{i}  \Delta v_{i} \otr u_{i};
\end{equation}

\begin{equation}
\label{eq:18}
\sum_{i}u_{i} \otr 1\otr v_{i}+\sum_{i} \sigma_{132}  (\Delta u_{i}\otr
v_{i} )
=\sum_{i} 1\otr v_{i}\otr u_{i}+\sum_{i} u_{i}\otr \Delta v_{i}.
\end{equation}

Combining Equations (\ref{eq:16})-(\ref{eq:18})
((\ref{eq:16})+(\ref{eq:17})-(\ref{eq:18})), we  obtain:

$$\sum_{i} (\Delta u_{i}\otr v_{i} -\Delta v_{i} \otr u_{i})=\sum_{i}
(1\otr  u_{i}\otr v_{i}-1  \otr v_{i}\otr u_{i} +u_{i}\otr 1\otr v_{i}
- v_{i} \otr  1\otr u_{i}), $$
which is equivalent to Equation  (\ref{eq:delta}).
Here we  used the identity:
$  \sigma_{231}( \Delta u_{i}\otr v_{i} ) =\sigma_{132}  (\Delta u_{i}\otr
v_{i}
   )$, which
can be easily verified  using the fact that
   $\Delta u_{i}$ is symmetric.

The final conclusion essentially follows from
Equation (\ref{eq:delta}). To see this, let us write
$\alt T=\sum_{i} u_{i}  \otr v_{i}$,
where $\{v_{i}\in UA\}$ are  assumed to be $R$-linearly independent.

  From  Equation (\ref{eq:delta}), it follows that
$$\sum_{i} \Delta u_{i}\otr v_{i}=
\sum_{i}  (1\otr u_{i}\otr v_{i}+u_{i}\otr 1\otr v_{i}). $$
I.e., $\sum_{i} ( \Delta u_{i}-1\otr u_{i}-u_{i}\otr 1)\otr v_{i}
=0$. Hence $\Delta u_{i}= 1\otr u_{i}+u_{i}\otr 1$, which
implies that $u_{i}\in \gm (A)$. Since $\alt T$ is
skew symmetric,  we conclude that $\alt T\in \gm (\wedge^{2} A)$. \qed
{\bf Remark. } It might  be useful to consider  the following cochain complex:
\begin{equation}
0\rightarrow
R\stackrel{\partial}{\rightarrow}UA\stackrel{\partial}{\rightarrow}
UA\otr UA\stackrel{\partial}{\rightarrow}UA\otr UA\otr UA
\stackrel{\partial}{\rightarrow}
\end{equation}
where $\partial :\otr^{n}UA\lon \otr^{n+1}UA, \ \partial=
\partial^{0}-\partial^{1}+\cdots +(-1)^{n+1}\partial^{n+1}, \ \
\partial^{i}(x_{1}\otr \cdots \otr x_{n})=\\
x_{1}\otr \cdots \otr x_{i-1}\otr \Delta x_{i}\otr x_{i+1} \otr \cdots \otr
x_{n}$
for $1\leq i\leq n$, and $\partial^{0}x=1\otr x, \ \partial^{n+1}x=x\otr 1$.
It is simple to check that $\partial^{2}=0$. In fact, this
is the subcomplex  of the  Hochschild cochain complex of
the algebra $C^{\infty}(G)$ ($G$ is a local Lie groupoid integrating the Lie
algebroid $A$) by restricting to  the space of
   left  invariant muti-differential
operators. It is natural to expect that the cohomology
of this complex is isomorphic to $\gm (\wedge^* A)$,
where the isomorphism from the cohomology group (more precisely, cocycles)
   to $\gm (\wedge^* A)$ is the usual skew-symmetrization  map.
   This is known to be true for Lie algebras \cite{dr:quasi} and
the tangent bundle Lie algebroid \cite{K}. However we could not find
such a    general result in the literature. In terms of
this cochain complex, it is simple to describe  what we have proved in
Proposition  \ref{pro:appendix}. It  simply means
   that $\alt : UA\otr UA \lon UA\otr UA$ maps
   2-cocycles into $\gm (\wedge^2 A)$.

We end this paper  by  a list of open questions.

{\bf Question 1}:   We  believe that   techniques in \cite{EK}
would be useful to prove the conjecture in Section 6.
While the proof of Etingof and Kazhdan relies heavily
on the double of a Lie bialgebra, the double of
a Lie bialgeboid is no longer a Lie algeboid.
Instead it is a Courant algebroid \cite{LWX}, where
certain anomalies are inevitable.  As a first step,
it is natural to ask: what is the universal enveloping algebra of a
Courant algebroid?
Roytenberg and Weinstein proved that Courant algebroids give
rise to homotopy Lie algebras \cite{RW}. It is expected that
these homotopy Lie algebras are useful to understand this
question as well as  the quantization problem.

{\bf Question 2}: One can form a    Kontsevich's formality-type conjecture
    for Lie algebroids, where one simply replaces in
Kontsevich formality theorem \cite{K} polyvector fields
by sections of $\wedge^*A$ and multi-differential operators by $UA\otr
   \cdots \otr UA$ for a Lie algebroid $A$.
   Does this conjecture hold?
It is not clear if the method in \cite{K} can be generalized
to the   context of general  Lie algebroids.
If this conjecture holds, it would imply that any
triangular Lie bialgebroid is quantizable.

{\bf Question 3}: Given a  solution  $r: \eta^* \lon
\frakg \ot \frakg$  of the classical dynamical Yang-Baxter
equation:
\begin{equation}
\alt dr -[r^{12} ,r^{13}]-[r^{12}, r^{23}]-[r^{13}, r^{23}]=0
\end{equation}
(in this case $\alt (r)$ satisfies Condition (ii) of a dynamical
$r$-matrix as   in  Section 2, if $r+r^{21}$ is ad-invariant),
a quantization of $r$ is a quantum dynamical $R$-matrix
$R: \eta^* \lon \uqg \ot \uqg $ such that $R(\lambda )=1 +\hbar {r}
(mod~\hbar^{2})$, where $\uqg$ is a quantum universal enveloping algebra.
   Is every classical  dynamical $r$-matrix quantizable? Many examples
are known
to be
quantizable (e.g.,  see \cite{ESS} for the
quantization of  classical dynamical $r$-matrices
in Schiffmann's classification list, and
    \cite{Xu4} for
the quantization of classical triangular dynamical $r$-matrices).
However, this problem still remains open for a general
dynamical $r$-matrix.

{\bf Question 4}:  According to the general principle  of
deformation theory, any deformation  corresponds to
a certain cohomology.
In particular, the deformation of a Hopf algebra is controlled
by  the cohomology of a certain  double complex
 \cite{GS1} \cite{GS2}. It is  natural
to ask what is the proper cohomology theory  controlling
the deformation of  a Hopf algebroid, and in particular
what is the premier obstruction to the quantization problem.

{\bf Question 5}:  What is the connection between
dynamical quantum groupoids and quantum Virasoro
algebra or quantum W-algebras \cite{FF} \cite{SKAO}?


\end{document}